\documentclass[12pt]{article}
\usepackage{amsmath,amsthm,amssymb,color,amscd}

\def\seq#1#2#3{#1_{#2},\,\ldots,#1_{#3}}

\def\w{\widetilde}

\def\b{\overline}

\def\PS{Poincar\'e series }

\def\tt{{\underline{t}}}

\def\mm{\underline{m}}

\def\1{\underline{1}}
\def\0{\underline{0}}
\def\R{\mathbb R}

\def\Q{\mathbb Q}
\def\N{\mathbb N}
\def\L{\mathbb L}

\def\Z{\mathbb Z}
\def\Q{\mathbb Q}
\def\C{\mathbb C}
\def\K{\mathbb K}

\def\OO{{\mathcal O}}
\def\EE{{\mathcal E}}

\def\CC{{\mathcal C}}

\def\calE{{\mathcal E}}

\def\DD{{\mathcal D}}

\def\oE{\stackrel{\circ}{E}}

\newtheorem{theorem}{Theorem}

\newtheorem{corollary}{Corollary}
\newtheorem{lemma}{Lemma}
\newtheorem{proposition}{Proposition}

\theoremstyle{definition}
\newtheorem{definition}{Definition}
\newtheorem{remark}{Remark}

\newtheorem{example}{Example}

\newenvironment{Proof}
{\noindent{\bf Proof\/}}{{ $\Box$}\smallskip\par}

\title{Poincar\'e series of valuations on functions over a subfield of complex numbers.
}

\author{
A.~Campillo,
\and F.~Delgado\thanks{The
first two authors partially supported by
grant PID2022-138906NB-C21 funded by
MICIU/AEI/10.13039/501100011033 and by ERDF/EU ``A way of making Europe''.}\; ,
\and S.M.~Gusein-Zade\thanks{
The work of the third author (sections \ref{sect:intro}, \ref{sect:resolution} and \ref{sect:PS-curve})
was supported by the grant 24-11-00124 of the Russian Science Foundation.
} }

\date{}
\begin{document}

\def\eps{\varepsilon}

\maketitle

\begin{abstract}
For a subfield $\K$ of the field $\C$ of complex numbers, we consider curve and
divisorial valuations on the algebra $\K[[x,y]]$ of formal power series in two
variables with the coefficients in $\K$. We compute the semigroup Poincar\'e series
and the classical Poincar\'e series of a (discrete, rank one)
valuation on it.
\end{abstract}

\section{Introduction}\label{sect:intro} 
A (discrete, rank one) valuation on the algebra
$\OO_{\C^2,0}$ of germs of holomorphic functions in two variables is a map
$\nu: \OO_{\C^2,0}\to \Z_{\ge 0}\cup \{+\infty\}$ such that
\begin{enumerate}
\item[1)] $\nu(0)=+\infty$;
\item[2)] $\nu(\lambda f)=\nu(f)$ for $f\in \OO_{\C^2,0}$,
$\lambda\in\C^*:=\C\setminus\{0\}$;
\item[3)] $\nu(f_1+f_2)\ge\min\{\nu(f_1),\nu(f_2)\}$ for $f_1,f_2\in \OO_{\C^2,0}$;
\item[4)] $\nu(f_1 f_2)=\nu(f_1)+\nu(f_2)$ for $f_1,f_2\in \OO_{\C^2,0}$.
\end{enumerate}
There are essentially two classes of valuations on
$\OO_{\C^2,0}$: curve valuations and divisorial
ones (see, e.\,g., \cite{Spiv}).

The notion of the Poincar\'e series of a collection $\{\nu_i\}$, $i=1,\ldots,r$,
of valuations was introduced in \cite{CDK}. For $r=1$ it reduces to
\begin{equation}\label{eq1}
P_{\nu}(t)=\sum_{v=0}^{\infty} \dim_{\C}(J^{\C} (v)/J^{\C} (v+1))\cdot t^v\in \Z[[t]],
\end{equation}
where
$J^{\C}(v):=\{f\in \OO_{\C^2,0} : \nu(f)\ge v\}$.
The Poincar\'e series of one curve valuation was known to specialists for many
years. Computations can be found, e.\,g., in~\cite{FAA-1999}.
The Poincar\'e series of one divisorial valuation was
computed in~\cite{Galindo-1995}.
The Poincar\'e series of a collection of curve valuations was computed in
\cite{Duke} (see also \cite{IJM} for a ``simpler'' proof).
It was shown that (for a collection consisting of more than one valuation) this
series coincides with the Alexander polynomial in several
variables of the corresponding algebraic link. This implies that the \PS determines
the topology of the union of the corresponding curves.

In~\cite{BLMS}, there were considered valuations on the algebra
$\EE_{\R^2,0}$ of germs of real analytic functions on $(\R^2,0)$. (All valuations on
$\EE_{\R^2,0}$ are restrictions of valuations on $\OO_{\C^2,0}$.)
There were considered several versions of the Poincar\'e series (including the one
defined in~\cite{CDK}). There was computed the Poincar\'e series for one curve valuation. In~\cite{MMJ},
there was computed the Poincar\'e series of one divisorial
valuation in this setting.

Here we consider Poincar\'e series of valuations on the algebra
$\EE_{\K^2,0}$ of germs of  holomorphic functions in two variables whose Taylor
coefficients are from a subfield $\K$ of the field
$\C$ of complex numbers (or, equivalently, on the algebra $\K[[x,y]]$
of formal power series). We compute the Poincar\'e series of one
curve or one divisorial valuation on it:
Theorems~\ref{theo:theo1}, \ref{theo:theo2}, and~\ref{theo:theo3}.
The Poincar\'e series are expressed in terms of a resolution
of a special form and are represented as products/ratios
of binomials of the form $(1-t^m)$.

The computation of the Poincar\'e series of a curve valuation
presented in~\cite{FAA-1999} was combinatorial. However,
it was known long before that. The importance of that paper
was the observation that it coincided with the monodromy
zeta function, i.\,e., with the Alexander polynomial of the
corresponding algebraic knot divided by $(1-t)$.
All the computations of the Poincar\'e series made by the authors
(starting from~\cite{Duke}) used some elements of topology:
versions of the integration with respect to the Euler characteristic.
The key property for that was the fact that the (additive)
Euler characteristic of a complex affine space is equal
to $1$. This does not work already for real analogues of the
Poincar\'e series in~\cite{BLMS} and~\cite{MMJ} (the Euler
characteristic of a real affine space is equal to $\pm 1$)
and makes no (straightfoward) sense for affine spaces over
subfields of $\C$. Therefore the computations in these cases
(at least at the contemporary stage) are essentially combinatorial
(modulo a certain analytic hint explained in Section~\ref{sect:PS-curve}).

The content of the paper is the following. In Section~\ref{sect:valuations}, we
discuss curve valuations
on the algebra $\EE_{\K^2,0}$. Section~\ref{sect:resolution} is
devoted to a description of a resolution process of an irreducible
curve playing the central role in the computations.
In Section~\ref{sect:subgroups}, one associates Galois groups
and subfields to the components of the exceptional divisor
of a resolution. Section~\ref{sect:semigroup} is devoted to the computation of
the semigroup Poincar\'e series, i.\,e., the sum
$\sum_{v\in S^{\K}_{\nu}}t^v$ over the set $S^{\K}_{\nu}$
of values of the valuation $\nu$ on $\EE_{\K^2,0}$: the semigroup
of values of the valuation. In Section~\ref{sect:PS-curve}, the
classical Poincar\'e series~(\ref{eq1}) of a curve valuation
is computed. In Section~\ref{sect:PS-divis} it is computed for
a divisorial valuation.

We are very thankful to the referee for reading carefully
the paper and making a number of useful suggestions.

\section{Curve valuations on functions over a field}\label{sect:valuations}
Let $\K$ be a subfield of the field $\C$ of complex numbers
and let $\EE_{\K^2,0}$ be the algebra of germs of holomorphic functions in two variables with the Taylor coefficients from $\K$.
A (discrete,
rank one) valuation on $\EE_{\K^2,0}$ (defined as above) is the restriction of a
valuation on $\OO_{\C^2,0}$: see, e.\,g., \cite{END}.
In the framework of the discussion below,
without changes it is possible to consider,
instead of the algebra $\EE_{\K^2,0}$, the algebra
$\K[[x,y]]$ of formal power series in two variables.
We shall refer to $\EE_{\K^2,0}$ because of some convenience
of the language.

\begin{remark}
Strictly speaking, some constructions below (in particular the description
of the action of the Galois group $G$ on branches at the beginning of
section~\ref{sect:resolution}) cannot be applied directly to the algebra
$\EE_{\K^2,0}$. The series for the branch $g\gamma$ and for the function $gh$ therein
could be non-convergent. These constructions work in the algebra
$K[[x,y]]$, however use of it everywhere produces some language problems. We do not
pay attention to this problem since, for a bounded values of valuations we can treat
not functions or series, but jets of order high enough, what means essentially
polynomials. Informally one can say that the computations can be made in a sort of
``an algebra of jets".
\end{remark}

Let $C$ be a plane algebroid branch on
$(\C^2,0)$, that is $C$ (after an appropriate change of the coordinates in $\C^2$ defined over $\K$) is given by
\begin{equation}\label{eq2}
x=x(\tau)=\tau^m,\\
y=y(\tau)=\sum_{i\ge m} c_i\tau^i\in \C[[\tau]]\,
\end{equation}
with some $m>0$.
(We assume that the parametrization (\ref{eq2}) is primitive, i.\,e., the
greatest
common divisor
of $m$ and of all
$i$ such that $c_i\ne 0$ is equal to $1$.) The branch $C$ defines a valuation
$\nu=\nu_C$ on the algebra $\OO_{\C^2,0}$ of holomorphic function germs in $x$ and $y$.  For
$f\in\OO_{\C^2,0}$, $\nu(f)$ is the degree of the leading term in the series
$f(x(\tau), y(\tau))\in\C[[\tau]]$:
\begin{equation}\label{eqn:curve_valuation}
 f(x(\tau), y(\tau))=a(f)\cdot\tau^{\nu(f)}+
 \text{ terms of higher degree}\;,
\end{equation}
where $a(f)\ne 0$.
If $f(x(\tau), y(\tau))\equiv 0$, $\nu(f)):=+\infty$.
Let us consider its restriction to the subalgebra
$\EE_{\K^2,0}\subset\OO_{\C^2,0}$. For $v\in\Z_{\ge 0}$, let
$J^\K(v):=\{f\in \EE_{\K^2,0}: \nu(f)\ge v\}$.
If all the quotients  $J^{\K}(v)/J^{\K}(v+1)$ are finite-dimensional (over $\K$) the
definition in \cite{CDK} reduces to the following one:

\begin{definition}
 The (classical) $\K$-{\em Poincar\'e series} of the branch $C$ is
 \begin{equation}\label{eqn:K-Poincare}
P_C^{\K}(t) =\sum_{v=0}^{\infty}\dim_{\K}\left(J^{\K}(v)/J^{\K}(v+1)\right)\cdot
t^v.
 \end{equation}
\end{definition}

 If $\K=\R$, $P^{\R}_C(t)$ is the {\em classical Poincar\'e series}
 of the valuation $\nu$ computed in~\cite{BLMS}.

\begin{remark}
 If $\K=\C$, all the coefficients of the
 series~(\ref{eqn:K-Poincare}) are equal to $0$ or $1$.
 This is not the case, in general, for a proper subfield
 $\K\subset\C$. In particular, if
 $\K=\R$ (the field of real numbers),
 all the coefficients are equal to $0$, $1$, or $2$
 (\cite{BLMS}).
\end{remark}

Let $S^\K_{C}:= \{\nu_C(f) : f\in \EE_{\K^2,0}\} \subset \Z_{\ge 0}$ be the {\em semigroup
of values} of the valuation $\nu_C$ on $\EE_{\K^2,0}$.

\begin{definition}\label{defsmg}
The $\K$-{\em semigroup Poincar\'e series} of the valuation $\nu$ is
$P^{S}_{C}(t) = \sum_{v\in S^{\K}_C} t^v$.
\end{definition}

 For $\K=\R$ this is the semigroup Poincar\'e series considered (and computed) in \cite{BLMS}.

\begin{example} (cf.~\cite[Example~1 and Remark~3]{BLMS})
 Let $\K=\R$, and let the branch $C$ be given by the
parametrization $(\tau^4, \alpha \tau^6+\tau^7)$ with a generic complex
$\alpha$. The semigroup of values of the valuation $\nu_C$
on $\EE_{\R^2,0}$ is equal to $S_C^{\R}=\langle 4, 6, 25\rangle$.
Here $\langle \bullet, \bullet, \ldots\rangle$ denotes the subsemigroup of
$\Z_{\ge0}$ generated by the corresponding natural numbers:
$$
\langle m_1, \ldots, m_p\rangle=\{m\in \Z_{\ge0}:
m=k_1m_1+\ldots+k_pm_p \text{\ \ with\ \ } k_i\in\Z_{\ge0}\}.
$$
This semigroup coincides with the semigroup of values of the
curve valuation $\nu_{C'}$ on $\OO_{\C^2,0}$, where the curve
$C'$ is given by the parametrization $(\tau^4, \tau^6+\tau^{19})$.
A direct relation between the curves $C$ and $C'$ is not clear.
The semigroup Poincar\'e series (both of the valuation $\nu_C$ on
$\EE_{\R^2,0}$ and of the valuation $\nu_{C'}$ on
$\OO_{\C^2,0}$)
is given by
$$
P^S_C(t)=\frac{(1-t^{12})(1-t^{50})}{(1-t^4)(1-t^6)(1-t^{25})}\,.
$$
The classical Poincar\'e series of the valuation $\nu_C$ on
$\EE_{\R^2,0}$ is given by
$$
P^{\R}_C(t)=\frac{(1-t^{24})(1-t^{50})}{(1-t^4)(1-t^6)(1-t^{25})}\,.
$$
 \end{example}

\section{The resolution process of a branch}\label{sect:resolution}  
Let $G$ be the Galois group of the extension $\C/\K$. i.\,e., the group of
automorphisms of the field $\C$ which are trivial on $\K$.
The group $G$ acts on branches (irreducible complex plane curve germs).
This action can be defined by any of the following three ways.
\begin{enumerate}
 \item[1)] Let a branch $\gamma$ be given by a parametrization
 $x=\tau^k$, $y=\sum\limits_{i\ge k} a_i\tau^i$. Then the branch $g\gamma$, $g\in G$,
is given by the parametrization
 $x=\tau^k$, $y=\sum\limits_{i\ge k} g(a_i)\tau^i$.
 \item[2)] Let a branch $\gamma$ be given by an equation
 $h(x,y)=0$, where $h(x,y)=\sum\limits_{i,j\ge 0} b_{ij}x^iy^j$.
 Then the branch $g\gamma$ is given by the equation $gh(x,y)=0$,
 where $gh(x,y):=\sum\limits_{i,j\ge 0} g(b_{ij})x^iy^j$.
 \item[3)] Let  a branch $\gamma$ be considered as a subset
 of $(\C^2,0)$: $(\gamma,0)\subset (\C^2,0)$. The group $G$
 acts on the plane $\C^2$ by $g(x,y)=(g(x),g(y))$.
 Then the branch $g\gamma$ as a subset of $(\C^2,0)$ is just
 the image of $\gamma$ under this action.
\end{enumerate}

\begin{remark}
 For fixed $x$ and $y$ the parametrization of $\gamma$ is not
 unique. One can obtain another parametrization from one of them
 multiplying the parameter $\tau$ by a root $\xi$ of degree $k$ of $1$.
 For $g\in G$, $g(\xi)$ is also a root of degree $k$ of $1$.
 Thus this only leads to another parametrization the
 branch $g\gamma$. It is not difficult to see that, for a
 $G$-invariant curve germ $\CC$ (possibly reducible), there
 exists a function germ $h\in \EE_{\K^2,0}$ such that
 $\CC=\{h=0\}$.
\end{remark}

Let us consider the minimal $G$-invariant resolution process of the branch $C$ defining the valuation
(by a sequence of blow-ups). This means the following.
When we blow-up the origin in $\C^2$, we have got one of the following three possibilities.
\begin{enumerate}
 \item[1)] There are finitely many different branches
 in the $G$-orbit of the branch $C$, their strict transforms
 are smooth and intersect the exceptional
 divisor transversally at pairwise different points. We get the (minimal embedded)
$G$-invariant resolution of $C$ and the process is completed.
 \item[2)] The strict transforms of the branches from the $G$-orbit of $C$ (whose number may be infinite) intersect
 the exceptional divisor at
 finite number of points, but either they are not smooth, or do not intersect the
 exceptional divisor transversally, or some (different)
 of them intersect it at the same point.
(In this case, if the $G$-orbit of $C$ has infinitely many branches, the number of strict transforms intersecting the
exceptional divisor at one point is also infinite.)
In this case one has to continue the resolution process.
 \item[3)] The strict transforms of the branches from the $G$-orbit of
${C}$
 intersect the exceptional divisor at infinitely many different points.
 In this case the process is interrupted. (Do not mix with the
 term ``completed''. In this case, one cannot get a $G$-invariant
 resolution since blowing-ups infinitely many points
 does not produce a complex surface.)
\end{enumerate}

In Case 2) one has to continue the resolution process
blowing-up (all) the intersection points of the strict transforms
of the branches from the $G$-orbit of $C$ with the exceptional divisor (their
number being finite). One gets one of the analogues of the cases 1)--3) again. Etc.

As the result, we arrive to one of the following three situations.
\begin{enumerate}
 \item[I.] There are finitely many branches in the $G$-orbit of $C$
 and they are resolved at a certain step.
 \item[II.] The $G$-orbit of $C$ has infinitely many branches, but
 the described process goes on without interruption.
 This means that at each step one has to blow-up finitely many
 points, but there are infinitely many steps like that.
 One does not arrive to a $G$-invariant resolution of $C$,
 but has got a $G$-invariant
resolution process (an infinite one).
 \item[III.] The process is interrupted at a certain step. This means
 that at this step the strict transforms of the branches from the $G$-orbit of $C$ intersect
 each of the new-born components of the exceptional divisor at
 infinitely many points.
 In this case one gets neither a $G$-invariant resolution of $C$,
 nor a resolution process.
\end{enumerate}

The cases I--III take place in the following situations.
\begin{enumerate}
\item[I.] All the coefficients in the series (\ref{eq2}) are from a finite extension
of $\K$. This is equivalent to the condition that the $G$-orbit of $C$ has finitely many elements.
\item[II.] The coefficients all together are not from a finite extension of $\K$, but
each of them is algebraic over $\K$.
\item[III.] There exists (at least) one coefficient transcendental over $\K$.
\end{enumerate}

At each step of a $G$-invariant resolution of $C$, the action
of the Galois group $G$ on $\C^2$ lifts to its action on the
surface of modification and thus induces an action on the set
of irreducible components of the exceptional divisor.
These components are in bijection with the vertices of
the dual graph of the modification and thus one has an action
of $G$ on this dual graph.

In the case I, the dual graph $\Gamma$ of the resolution looks
like in Figure~\ref{fig1}.
\begin{figure}[h]
$$
\unitlength=0.50mm
\begin{picture}(120.00,110.00)(0,-30)
\thinlines
\put(-30,30){\circle*{2}}
\put(-40,24){{\scriptsize ${\bf 1}=\sigma_0$}}

\put(-30,30){\line(1,0){90}}
\put(60,30){\circle*{2}}

\put(-15,30){\line(0,-1){15}}
\put(-15,30){\circle*{2}}
\put(-15,15){\circle*{2}}
\put(-16.5,10){{\scriptsize$\sigma_1$}}
\put(-16,33){{\scriptsize$\tau_1$}}

\put(5,30){\line(0,-1){20}}
\put(5,30){\circle*{2}}
\put(5,10){\circle*{2}}
\put(6.5,7){{\scriptsize$\sigma_2$}}
\put(4,33){{\scriptsize$\tau_2$}}

\put(35,30){\line(0,-1){20}}
\put(35,30){\circle*{2}}
\put(35,10){\circle*{2}}
\put(36.5,7){{\scriptsize$\sigma_3$}}
\put(34,33){{\scriptsize$\tau_3$}}

\put(50,33){\scriptsize{$\rho_1$}}

\put(60,30){\line(1,1){40}}
\put(85,55){\line(1,-1){10}}
\put(85,55){\circle*{2}}
\put(100,70){\circle*{2}}

\put(100,70){\line(2,1){20}}
\put(120,80){\vector(1,0){10}}
\put(120,80){\circle*{2}}
\put(110,75){\line(1,-2){4}}
\put(120,80){\line(1,-2){4}}

\put(100,70){\line(2,-1){20}}
\put(120,60){\vector(1,0){10}}
\put(120,60){\circle*{2}}
\put(110,65){\line(-1,-2){4}}
\put(120,60){\line(-1,-2){4}}

\put(60,30){\line(1,-1){40}}
\put(100,-10){\circle*{2}}
\put(85,5){\line(-1,-1){10}}
\put(85,5){\circle*{2}}
\put(86,6){{\scriptsize$\tau_4$}}
\put(75,-5){\circle*{2}}
\put(73,-13){{\scriptsize$\sigma_4$}}

\put(98,-5){{\scriptsize$\rho_2$}}

\put(100,-10){\line(2,1){20}}
\put(120,0){\vector(1,0){10}}
\put(120,0){\circle*{2}}
\put(110,-5){\line(1,-2){4}}
\put(120,0){\line(1,-2){4}}

\put(100,-10){\line(2,-1){20}}
\put(120,-20){\vector(1,0){10}}
\put(120,-20){\circle*{2}}
\put(120,-20){\line(-1,-2){4}}
\put(130,-17){{\scriptsize $C$}}
\put(110,-15){\line(-1,-2){4}}
\put(116,-28){\circle*{2}}
\put(118,-15){{\scriptsize$\tau_6$}}
\put(110,-15){\circle*{2}}
\put(106,-23){\circle*{2}}
\put(98,-28){{\scriptsize$\sigma_5$}}
\put(108,-31){{\scriptsize$\sigma_6$}}

\put(60,30){\line(1,0){40}}
\put(85,30){\line(0,-1){13}}
\put(85,30){\circle*{2}}
\put(100,30){\circle*{2}}

\put(100,30){\line(2,1){20}}
\put(120,40){\line(1,-2){4}}
\put(120,40){\vector(1,0){10}}
\put(120,40){\circle*{2}}

\put(110,35){\line(1,-2){4}}
\put(100,30){\line(2,-1){20}}
\put(120,20){\line(-1,-2){4}}
\put(120,20){\vector(1,0){10}}
\put(120,20){\circle*{2}}
\put(110,25){\line(-1,-2){4}}

\end{picture}
$$
\caption{The resolution graph $\Gamma$ in case I (here $s=2$, $g=6$, see below).}
\label{fig1}
\end{figure}

For a description of the resolution (or of the resolution process),
it is convenient to use the quotient
$\check \Gamma$ of the graph $\Gamma$ by the action of the
Galois group $G$. In this case it looks like in Figure~\ref{fig2}.

The graph $\Gamma$ contains a subgraph which geometrically coincides with the dual graph of a usual
resolution of the curve $C$: the lower part of $\Gamma$ on Figure \ref{fig1}. (We do not consider this resolution below.) We use
``short" notations (by greek letters) for vertices on this subgraph (and also for
the corresponding vertices of the graph $\check \Gamma$). All other vertices are
obtained from them by the action of elements of the Galois group $G$.

In Figures \ref{fig1} and \ref{fig2},
$\sigma_i$, $i=0,1,\ldots, g$, (and also their images under the action of $G$) are
the {\it dead ends\/} of the graph ($g$ is the number of the Puiseux pairs of
the curve $C$),
$\tau_i$,
$i=1,\ldots,g$, (and also their images) are the {\it rupture\/} points of it,
$\rho_j$, $j=1,\ldots,s$,
are the {\it
splitting points\/} corresponding to the changes of the isotropy subgroups
for the $G$-action on $\Gamma$ (that is the splitting points for usual
resolutions of
different branches from the orbit of $C$, there are finitely many of them in Case I).
Let $\delta_C$ be the vertex corresponding to the
component of the exceptional divisor intersecting the strict
transform of the curve $C$.
In Figure~\ref{fig1}, $\delta_C=\tau_6$.
Figure~\ref{fig2} shows $\check \Gamma$ in the case when the process of the
$G$-resolution of $C$ is not finished when $C$ itself is resolved, $\delta_C=\rho_s$
there.
(Otherwise $\delta_C$ coincides with $\tau_g$.)

\begin{figure}[h]
$$
\unitlength=0.50mm
\begin{picture}(120.00,40.00)(0,10)

\thinlines
\put(-30,30){\circle*{2}}
\put(-40,24){{\scriptsize ${\bf 1}=\sigma_0$}}

\put(-30,30){\line(1,0){20}}
\put(-8,30){\circle*{0.5}}
\put(-5,30){\circle*{0.5}}
\put(-2,30){\circle*{0.5}}
\put(1,30){\circle*{0.5}}
\put(4,30){\circle*{0.5}}
\put(5,30){\line(1,0){40}}

\put(-15,30){\line(0,-1){15}}
\put(-15,30){\circle*{2}}
\put(-15,15){\circle*{2}}
\put(-16.5,10){{\scriptsize$\sigma_1$}}
\put(-16,33){{\scriptsize$\tau_1$}}

\put(10,30){\line(0,-1){20}}
\put(10,30){\circle*{2}}
\put(10,10){\circle*{2}}
\put(11.5,7){{\scriptsize$\sigma_q$}}
\put(9,33){{\scriptsize$\tau_q$}}

\put(30,30){\circle*{2}}
\put(27,23){\scriptsize{$\rho_1$}}

\put(30,30){\line(1,1){10}}
\put(30,30){\line(1,2){8}}
\put(30,30){\line(0,1){10}}

\put(40,30){\line(0,-1){20}}
\put(40,30){\circle*{2}}
\put(40,10){\circle*{2}}
\put(41.5,7){{\scriptsize$\sigma_{q+1}$}}
\put(39,33){{\scriptsize$\tau_{q+1}$}}

\put(48,30){\circle*{0.5}}
\put(51,30){\circle*{0.5}}
\put(54,30){\circle*{0.5}}

\put(55,30){\line(1,0){40}}

\put(60,30){\line(0,-1){20}}
\put(60,30){\circle*{2}}
\put(60,10){\circle*{2}}

\put(75,30){\line(2,1){10}}
\put(75,30){\line(1,1){10}}
\put(75,30){\line(1,2){8}}
\put(75,30){\line(0,1){10}}

\put(75,30){\circle*{2}}
\put(72,23){\scriptsize{$\rho_2$}}

\put(90,30){\line(0,-1){20}}
\put(90,30){\circle*{2}}
\put(90,10){\circle*{2}}

\put(98,30){\circle*{0.5}}
\put(101,30){\circle*{0.5}}
\put(104,30){\circle*{0.5}}

\put(106,30){\line(1,0){25}}
\put(120,30){\circle*{2}}
\put(118,23){\scriptsize{$\rho_i$}}

\put(110,30){\line(0,-1){15}}
\put(110,30){\circle*{2}}
\put(110,15){\circle*{2}}
\put(111.5,12){\scriptsize{$\sigma_g$}}
\put(108,33){\scriptsize{$\tau_g$}}

\put(120,30){\line(2,1){10}}
\put(120,30){\line(1,1){10}}
\put(120,30){\line(0,1){10}}

\put(134,30){\circle*{0.5}}
\put(137,30){\circle*{0.5}}
\put(140,30){\circle*{0.5}}

\put(142,30){\line(1,0){8}}

\put(150,30){\vector(1,-1){10}}
\put(150,30){\circle*{2}}
\put(162,25){{\scriptsize $C$}}
\put(132,23){\scriptsize{$\rho_s = \delta_C$}}

\put(150,30){\line(2,1){10}}
\put(150,30){\line(1,1){10}}
\put(150,30){\line(1,2){8}}
\put(150,30){\line(0,1){10}}

\end{picture}
$$
\caption{The quotient $\check \Gamma$ of the resolution graph in Case I (with
arbitrary $g$ and $s$).}
\label{fig2}
\end{figure}

In Case II, one has the dual graph of the resolution process.
The corresponding quotient $\check \Gamma$
looks like in Figure~\ref{fig3}. One can say that $s$ is equal to infinity in this
case.

\begin{figure}[h]
$$
\unitlength=0.50mm
\begin{picture}(120.00,40.00)(0,10)

\thinlines
\put(-30,30){\circle*{2}}
\put(-40,24){{\scriptsize ${\bf 1}=\sigma_0$}}

\put(-30,30){\line(1,0){20}}
\put(-8,30){\circle*{0.5}}
\put(-5,30){\circle*{0.5}}
\put(-2,30){\circle*{0.5}}
\put(1,30){\circle*{0.5}}
\put(4,30){\circle*{0.5}}
\put(5,30){\line(1,0){47}}

\put(-15,30){\line(0,-1){15}}
\put(-15,30){\circle*{2}}
\put(-15,15){\circle*{2}}
\put(-16.5,10){{\scriptsize$\sigma_1$}}
\put(-16,33){{\scriptsize$\tau_1$}}

\put(10,30){\line(0,-1){20}}
\put(10,30){\circle*{2}}
\put(10,10){\circle*{2}}
\put(11.5,7){{\scriptsize$\sigma_q$}}
\put(9,33){{\scriptsize$\tau_q$}}

is

\put(30,30){\circle*{2}}
\put(27,23){\scriptsize{$\rho_1$}}

\put(30,30){\line(1,1){10}}
\put(30,30){\line(1,2){8}}
\put(30,30){\line(0,1){10}}

\put(40,30){\line(0,-1){20}}
\put(40,30){\circle*{2}}
\put(40,10){\circle*{2}}
\put(41.5,7){{\scriptsize$\sigma_{q+1}$}}
\put(39,33){{\scriptsize$\tau_{q+1}$}}

\put(48,30){\circle*{0.5}}
\put(51,30){\circle*{0.5}}
\put(54,30){\circle*{0.5}}

\put(55,30){\line(1,0){43}}

\put(60,30){\line(0,-1){20}}
\put(60,30){\circle*{2}}
\put(60,10){\circle*{2}}

\put(75,30){\circle*{2}}
\put(72,23){\scriptsize{$\rho_j$}}

\put(75,30){\line(2,1){10}}
\put(75,30){\line(1,1){10}}
\put(75,30){\line(1,2){8}}
\put(75,30){\line(0,1){10}}

\put(90,30){\line(0,-1){20}}
\put(90,30){\circle*{2}}
\put(90,10){\circle*{2}}
\put(91.5,7){{\scriptsize$\sigma_{g}$}}
\put(89,33){{\scriptsize$\tau_{g}$}}

\put(98,30){\circle*{0.5}}
\put(101,30){\circle*{0.5}}
\put(104,30){\circle*{0.5}}
\put(107,30){\circle*{0.5}}

\put(112,30){\line(1,0){20}}
\put(120,30){\circle*{2}}

\put(120,30){\line(2,1){10}}
\put(120,30){\line(1,1){10}}
\put(120,30){\line(0,1){10}}

\put(134,30){\circle*{0.5}}
\put(137,30){\circle*{0.5}}
\put(140,30){\circle*{0.5}}

\put(142,30){\line(1,0){12}}
\put(150,30){\circle*{2}}
\put(147,23){\scriptsize{$\rho_k$}}

\put(150,30){\line(2,1){10}}
\put(150,30){\line(1,1){10}}
\put(150,30){\line(1,2){8}}
\put(150,30){\line(0,1){10}}

\put(156,30){\circle*{0.5}}
\put(159,30){\circle*{0.5}}
\put(162,30){\circle*{0.5}}
\put(165,30){\circle*{0.5}}

\end{picture}
$$
\caption{The quotient $\check \Gamma$ of the graph of the resolution process in Case II (with $s=\infty$).}
\label{fig3}
\end{figure}

In Case III, the quotient of the modification graph
looks like in Figure~\ref{fig4}.
The curve $C$ may be not resolved at the interruption step.

In the latter case, the valuation coincides with $n\nu$, where $\nu=\nu_{\rho_s}$ is the
divisorial valuation defined by $E_{\rho_s}$ and $n$ is the intersection multiplicity
$E_{\rho_s}\circ \w C$ of the strict transform $\w C$ of the curve $C$ with
the component $E_{\rho_s}$.


\begin{figure}[h]
$$
\unitlength=0.50mm
\begin{picture}(120.00,40.00)(-20,0)

\thinlines
\put(-30,30){\circle*{2}}
\put(-40,24){{\scriptsize ${\bf 1}=\sigma_0$}}

\put(-30,30){\line(1,0){20}}
\put(-8,30){\circle*{0.5}}
\put(-5,30){\circle*{0.5}}
\put(-2,30){\circle*{0.5}}
\put(1,30){\circle*{0.5}}
\put(4,30){\circle*{0.5}}
\put(5,30){\line(1,0){40}}

\put(-15,30){\line(0,-1){15}}
\put(-15,30){\circle*{2}}
\put(-15,15){\circle*{2}}
\put(-16.5,10){{\scriptsize$\sigma_1$}}
\put(-16,33){{\scriptsize$\tau_1$}}

\put(10,30){\line(0,-1){20}}
\put(10,30){\circle*{2}}
\put(10,10){\circle*{2}}
\put(11.5,7){{\scriptsize$\sigma_q$}}
\put(9,33){{\scriptsize$\tau_q$}}

\put(30,30){\circle*{2}}
\put(27,23){\scriptsize{$\rho_1$}}

\put(30,30){\line(1,1){10}}
\put(30,30){\line(1,2){8}}
\put(30,30){\line(0,1){10}}

\put(40,30){\line(0,-1){20}}
\put(40,30){\circle*{2}}
\put(40,10){\circle*{2}}
\put(41.5,7){{\scriptsize$\sigma_{q+1}$}}
\put(39,33){{\scriptsize$\tau_{q+1}$}}

\put(48,30){\circle*{0.5}}
\put(51,30){\circle*{0.5}}
\put(54,30){\circle*{0.5}}

\put(55,30){\line(1,0){30}}

\put(60,30){\line(0,-1){20}}
\put(60,30){\circle*{2}}
\put(60,10){\circle*{2}}

\put(75,30){\circle*{2}}
\put(72,23){\scriptsize{$\rho_s$}}

\put(75,30){\line(2,1){10}}
\put(75,30){\line(1,1){10}}
\put(75,30){\line(1,2){8}}
\put(75,30){\line(0,1){10}}
\put(75,30){\line(2,-1){8}}

\put(90,30){{\scriptsize{\bf $\infty$}}}

\put(77,28){\circle*{0.5}}
\put(79,26){\circle*{0.5}}
\put(81,24){\circle*{0.5}}
\put(83,22){\circle*{0.5}}
\put(85,20){\circle*{0.5}}
\put(87,18){\circle*{0.5}}
\put(89,16){\circle*{0.5}}

\put(90,15){\circle*{0.5}}
\put(90,12){\circle*{0.5}}
\put(90,10){\circle*{0.5}}
\put(90,8){\circle*{0.5}}
\put(90,5){\circle*{0.5}}
\put(90,2){\circle*{0.5}}

\put(93,15){\circle*{0.5}}
\put(96,15){\circle*{0.5}}
\put(99,15){\circle*{0.5}}
\put(102,15){\circle*{0.5}}
\put(105,15){\circle*{0.5}}
\put(108,15){\circle*{0.5}}
\put(110,15){\circle*{0.5}}
\put(113,15){\circle*{0.5}}
\put(115,15){\circle*{0.5}}

\put(115,15){\circle*{0.5}}
\put(115,12){\circle*{0.5}}
\put(115,10){\circle*{0.5}}
\put(115,8){\circle*{0.5}}
\put(115,5){\circle*{0.5}}
\put(115,2){\circle*{0.5}}

\put(117,17){\circle*{0.5}}
\put(119,19){\circle*{0.5}}
\put(121,21){\circle*{0.5}}
\put(121,21){\vector(1,1){5}}

\end{picture}
$$
\caption{The quotient $\check \Gamma$ of the modification graph in Case III.}
\label{fig4}
\end{figure}

\begin{example}
Let $\K=\Q$ and the branch $C$ given by the parametrization
$x= \tau^4$, $y = \alpha_0 \tau^6 + t^7 + \sum_{i\ge 1}\alpha_i \tau^{7+i}$ where
$\alpha_0, \alpha_1, \ldots \in \C\setminus \{0\}$.
Let $\K' :=\Q(\alpha_0, \alpha_1, \ldots)\subset \C$.

\begin{figure}[h]
$$
\unitlength=1.00mm
\begin{picture}(120.00,20.00)(-10,5)

\thinlines

\put(10,25){\line(1,0){80}}

\put(10,25){\circle*{2}}
\put(5,20){{${\bf 1}$}}
\put(30,25){\circle*{2}}
\put(30,25){\line(0,-1){20}}
\put(30,5){\circle*{2}}
\put(25,20){{${\bf 3}$}}
\put(25,5){{${\bf 2}$}}

\put(50,25){\circle*{2}}

\put(50,25){\line(0,-1){20}}
\put(50,5){\circle*{2}}
\put(45,20){{${\bf 5}$}}
\put(45,5){{${\bf 4}$}}

\put(60,25){\circle*{2}}
\put(70,25){\circle*{2}}
\put(80,25){\circle*{2}}
\put(90,25){\circle*{2}}

\put(59,20){{${\bf 6}$}}
\put(69,20){{${\bf 7}$}}
\put(79,20){{${\bf 8}$}}
\put(89,20){{${\bf 9}$}}

\put(92,25){\circle*{0.2}}
\put(94,25){\circle*{0.2}}
\put(96,25){\circle*{0.2}}
\put(98,25){\circle*{0.2}}
\put(100,25){\circle*{0.2}}
\put(102,25){\circle*{0.2}}
\put(104,25){\circle*{0.2}}

%
%
%
%

\end{picture}
$$
\caption{The dual graph of $C$}
\label{fig5}
\end{figure}

One can distinguish the following cases:

\begin{enumerate}
\item $\K'$ is finite algebraic extension of $\Q$. In this case
$\K' = \Q(\seq{\alpha}0k)$ for some $k$ and if we add the conditions
$\alpha_0\notin \Q$ (algebraic, generic)  and
$\alpha_i\notin \Q(\seq{\alpha}0{i-1})$ (algebraic, generic) for $1\le i \le k$ then
$s=k+1$ and the
 splitting points
are $\rho_1= 3$, $\rho_i = 5+(i-1)$ for $i=2,\ldots, k+1$.
Thus we are in the Case I.

Pay attention that, if $\alpha_0\in \K'$ is not generic, the pictire can be somewaht
different.
For example, let
$y = \sqrt{2} \tau^6 +\tau^7$. Then the branch $C$ has the equation
$f(x,y) = (y^2- 2 x^3)^2 -4\sqrt{2} y x^5 -x^7 = 0$ and the conjugate of $C$ with
respect to
$G$, $g C$ ($g\in G$ is given by $g(\sqrt{2})=-\sqrt{2}$),
has the parametrization $(\tau^4, -\sqrt{2}\tau^6 + \tau^7)$ and the equation
$f'(x,y) = g f(x,y) = (y^2- 2 x^3)^2 + 4\sqrt{2} y x^5 -x^7 = 0$. In this case the
only splitting point
is $\rho = 5$ (also $\delta_C=5$). Notice that the intersection multiplicity of $C$
and $C' = gC$ is $26$ and that the coordinate of the intersection point of the strict
transform of $C$
(resp. $C'$) is a certain $a\in \Q(\sqrt{2})\setminus \Q$ (resp. $g(a)$).

\item $\K'$ is an infinite algebraic extension of $\Q=\K$. Then $s=\infty$ and
we are in Case II. As in the previous case, in general, the splitting points are
$3, 6, 7, \ldots$.

\item $\K'$ is a trascendental extension of $\Q$. In this case we are in Case III.
For simplicity, let us assume that $\alpha_0$ is trascendental over $\Q$.
Then the process is interrumpted in the divisor $\delta_\nu=3$.
In order to get a function $g\in \OO_{\C^2,0}$ such that $\nu(g)=2 r +1 \ge
7$ one needs that
$$
g = (y^2- \alpha_{0}^2 x^3) + \text{ higher order terms } \; .
$$
But in this case the equation can not be rational. So, in
this case,
the semigroup $S^{\Q}_{\nu}$ is
$\langle 4, 6 \rangle \subset 2 \Z$ and the valuation
$\nu$ is
$2$ times the divisorial valuation  defined by the component $E_3$.
\end{enumerate}

\end{example}

\begin{remark}\label{rational}
Let $F\in \K[[x,y]]$ be irreducible.
The local ring of the irreducible curve (over $\K$) defined by $F$ is
$R = \K[[x,y]]/(F)$. The integral closure $\b{R}$ in its quotient field is the
valuation ring $R_\nu$
of a discrete valuation $\nu$, its residual field
$\K'$ is a finite
algebraic extension of $\K$ and could be realized and therefore identified with a
subfield of
$\C$.
The equation $F=0$ defines a reduced algebroid curve over the complex field $\C$, each
branch corresponds to an irreducible factor of $F$ in $\C[[x,y]]$. Let $C$ be one of
these branches. The valuation $\nu_C$ defined by $C$ on
$\K[[x,y]]$ coincides with the valuation $\nu$ (seen as a valuation on $\K[[x,y]]$
forgetting the quotient by $F$)
and the Poincar\'e series (over
$\K$) defined by $\nu_C$ is the Poincar\'e series defined by $\nu$. Note that all
the branches are equisingular, the minimal resolution of $\{F=0\}$ (in the complex
setting) coincides with the $G$-resolution of $\nu_C$ and all the valuations
$\nu_C$,
for complex branches $C$ of $\{F=0\}$, restrict to the same valuation $\omega$ in
$\K[[x,y]]$.
Moreover, if
$G'$ is the Galois group of the extension $\K\subset {\L}$, ${\L}$ the normal
closure of $\K'$ in $\C$, $G'$ is the quotient of $G$ by the isotropy group of the
field $\L$. This situation just corresponds to the case I described
above.

Let $C$ be a curve parametrized as in (\ref{eq2}). The map
$f\mapsto f(x(\tau), y(\tau))$ defines an homomorphism
$\Phi: \K[[x,y]]\to \C[[\tau]]$. If $\ker (\Phi)\neq 0$, then one has
$\ker(\Phi)=(F)$ for some irreducible $F\in \K[[x,y]]$ and we are in the above
described situation, i.\,e., also in Case I. Otherwise, $\ker(\Phi)=(0)$ and we are
in
Case II or III. In this case the valuation $\nu_C$ is centered in the two-dimensional
ring $\K[[x,y]]$. Case II  corresponds to the case in which the residual field
$\K_{\nu_C}$ of the valuation ring is an algebraic extension of $\K$ but not finite
dimensional. Case III corresponds to the case in which $\K_{\nu_C}$ is a
trascendental extension of $\K$.

\end{remark}

\section{Subgroups of the Galois group, subfields, divisors and curvettes}\label{sect:subgroups}
Further on we shall discuss only Cases I and II.
We do not discuss Case III since in this case the
valuation $\nu_C$ on $\EE_{\K^2,0}$ is a multiple
of a divisorial valuation whose Poincar\'e series
is described in Section~\ref{sect:PS-divis}.
Let $( E_\sigma \circ E_{\delta})$, $\sigma, \delta\in \Gamma$, be the
intersection matrix of the components $E_\sigma$ (infinite one in Case II). The
diagonal entry $E_\sigma\circ E_{\sigma}$ is the
self-intersection
number of the component $E_\sigma$ on the final modification in Case I and, in Case
II, a modification at a step, after which the points of $E_\sigma$ are not blown-up.
One has
$E_{\sigma}\circ E_{\sigma}<0$.
The non-diagonal entry $E_{\sigma}\circ E_{\delta}$ is equal to 1 (respectively to
$0$) if these components intersect (do not intersect respectively).

Let
$(m_{\sigma, \delta}) := - (E_\sigma \circ E_{\delta})^{-1}$.
In order not to treat infinite matrices,
in Case II the entry $m_{\sigma, \delta}$ can be computed
at a step of the resolution process after which the
components $E_{\sigma}$ and $E_{\delta}$ are not modified.

The number $m_{\sigma, \delta}$ can be interpreted as the intersection number of some
germs of curves on $(\C^2,0)$, so-called curvettes. A curvette corresponding to the
component $E_\sigma$ is the blow-down $C_\sigma$ of a complex-analytic smooth curve
$\gamma_\sigma$ (in Case II on a modification after which the component $E_\sigma$ is not modified) transversal to the
component $E_\sigma$ at a smooth point of the exceptional divisor, that is not an
intersection point with other components. One can show that $m_{\sigma, \delta}$ is
equal to the intersection number $C_\sigma \circ C_\delta$.

For $\sigma\in \Gamma$, let $m_\sigma=\nu_C(f_{\sigma})$,
where $f_\sigma=0$ is an equation defining a curvette
$C_{\sigma}$. It is equal to the intersection number
$C\circ C_\sigma$. One can see that, in Case I,
$m_\sigma$
is equal to $m_{\sigma, \delta_C}$ for $\delta_C$ such that, on the surface of the
resolution, the
strict transform of the branch $C$ intersects the component $E_{\delta_C}$.
In Case II, the vertex $\delta_C$ makes sense only at any
step of the resolution process of the curve $C$.
One can see that $m_\sigma=m_{\sigma, \delta}$
at a step of the resolution process
when the branch $C$ itself is resolved (not in the
$G$-invariant sense) and
after which the component $E_{\sigma}$ is not modified.

\subsubsection*{Isotropy subgroups and subfields.}
Let, as above, $G_\sigma$ be the
isotropy subgroup of the component $E_\sigma$, i.\,e., the isotropy subgroup of the
vertex $\sigma$ in $\Gamma$. The subgroup $G_\sigma$ has a finite index in the
Galois group $G$. Let $\K_\sigma\subset \C$ be the invariant subfield of the group
$G_\sigma$.
One can see that, for $g\in G$, $G_{g \sigma}=g G_{\sigma}g^{-1}$,
$\K_{g \sigma}= g \K_{\sigma}$. For all the vertices $\sigma\in \Gamma$ in between
the initial one $\sigma_0 =$\;{\bf 1} and $\rho_1$ ($\rho_1$ included) the subgroups
$G_{\sigma}$ and the subfields $\K_\sigma$ are the same:
$G_{\sigma}=G=G_0$, $\K_{\sigma}=\K=\K_0$. For all vertices
$\sigma\in \Gamma$ in between  the vertex $\rho_j$ and $\rho_{j+1}$
($\rho_j$ excluded, $\rho_{j+1}$ included)
the isotropy subgroups $G_{\sigma}$ and the corresponding subfields $\K_{\sigma}$ are the same:
$G_\sigma= G_j$, $\K_\sigma=\K_j$.
In Case I, for all the vertices $\sigma\in \Gamma$ between $\rho_s$ and $\delta_C$ one has
$G_{\sigma}=G_s$, $\K_\sigma=\K_s$. For vertices from $G$-orbits of the discussed
ones, the corresponding subgroups are conjugate to $G_{\sigma}$ and the subfields
are obtained from $\K_\sigma$ by the shift by an element of the Galois group $G$. One has
$G=G_0\supset G_1\supset \ldots \supset G_s$,
$\K=\K_0\subset \K_1\subset \ldots \subset \K_s$. One has $[\K_j: \K_{j-1}]=[G_{j-1}:G_j]$.

\subsubsection*{Standard charts.}
For computations, we shall use coordinate systems compatible with the actions of the
Galois groups. This means the following.
There are two standard charts covering the initial component of the resolution
process: with the coordinates $u=x, w=y/x$ or $u=y, w=x/y$ respectively.
If, on a previously created component one makes the blow-up at the point
$(u,w)=(0,a)$ (in a standard chart), one has two standard charts covering the
new-born component: with the coordinates $u'=u, w'=(w-a)/u$ and
$u'=w-a, w'=u/(w-a)$ respectively.

On all the components $E_\sigma$ of the exceptional divisor except $E_{\rho_j}$ and
those from its $G$-orbit, on each step the blow-up is made at one point
with the coordinates $(0, a)$ in a standard chart with
$a\in \K_{\sigma}$.
In this case it can be convenient to change the coordinate $w$
by $\widetilde{w}=w/(w-a)$ or
$\widetilde{w}=(w-a)/w$ so that the coordinates of the intersection points with the
components on the geodesic between $\sigma_0$ and
$\delta_C$ become equal to $0$ and $\infty$.
On the component $E_{\rho_j}$,
the blow-ups are made at a certain point $(0, a)$ with
$a\in \K_j\setminus \K_{j-1}$ (that is not from
$\K_{\rho_j}=\K_{j-1}$) and at the  points from its
$G_{j-1}$-orbit. In this case a coordinate change as above cannot be used since it is
not defined over $\K_{j-1}$.

\subsubsection*{Special curvettes.}
In what follows, we shall use curvettes of special type.
In a standard chart covering an affine part of $E_\sigma$, let the (smooth
and
transversal) curve $\gamma_\sigma$ be defined over the field $\K_\sigma$. In other
words, the curve $\gamma_\sigma$ is invariant under the action of $G_\sigma$.
In particular one can take the curve $\{u=\tau$, $w=c=\text{const}\}$ in a standard
chart with
$c\in\K_\sigma$ ($u=0$ is the equation of $E_\sigma$ in the chart).

\begin{definition}
A $\K_\sigma$-curvette at the component $E_\sigma$ is the blow-down $C_\sigma$ of a
curve $\gamma_\sigma$ described above.
\end{definition}

\begin{remark}
One can see that $C_\sigma$ has a parametrization like (\ref{eq2}) with all the
coefficients from
$\K_\sigma$.
A $\K_\sigma$-curvette at the component $E_\sigma$ can be
defined by an equation $f_\sigma=0$ with $f_\sigma\in \EE_{\K_\sigma^2,0}$.
\end{remark}

\begin{definition}
A $G_\sigma$-curvette at the component $E_\sigma$ is the union of the curves
$g C_\sigma$ ($C_\sigma$ is a $\K_\sigma$-curvette at $E_\sigma$) over
representatives of the $G_\sigma$-classes in $G$.
\end{definition}

One can show that a $G_\sigma$-curvette at the component $E_\sigma$ can be defined by
an equation $F_\sigma=0$ with
$F_\sigma\in \EE_{\K^2,0}$. In fact one can take
$F_\sigma=\prod_{[g]\in G/G_\sigma} (gf_\sigma)$, where the product is over
representatives of the $G_\sigma$-classes in $G$.

\begin{remark}
 Somewhat similar type of symmetries of
 the resolution graph was considered in~\cite{MNach}. In that case one had an action of a group on the graph, but
 not on the underlying space $\C^2$ or on the space of functions. The series computed in this paper and therein are different.
 In fact the setting studied here was one of the motivations for~\cite{MNach}.
\end{remark}

\section{The semigroup Poincar\'e series of a curve valuation}
\label{sect:semigroup}  
Let $\nu (=\nu_C)$ be the valuation on $\EE_{\K^2,0}$ defined by an irreducible
curve
$C\subset
(\C^2,0)$ and assume that $\nu_C$ is a valuation of the type I or II described in Section
\ref{sect:resolution}.
Let $P^S_C(t) = \sum_{m\in S^\K_\nu} t^m$ be the
semigroup Poincar\'e series of $\nu_C$
(see Definition \ref{defsmg}).

For $\sigma\in \check \Gamma$, let $M_\sigma = \nu(F_{\sigma})$
where $F_\sigma=0$ is the equation defining a $G_\sigma$-curvette at the component
$E_\sigma$.
One can see that $M_\sigma=\sum_{[g]\in G/G_{\sigma}} m_{g \sigma}$. (Compare with
the definition of $M_\sigma$ in \cite{BLMS}, where $\K=\R$.)

\medskip
The aim of this section is to prove the following statement.

\begin{theorem}\label{theo:theo1}
$$
P^S_{C}(t) = \frac{\prod_{i=1}^g (1-t^{M_{\tau_i}})}
{\prod_{i=0}^g (1-t^{M_{\sigma_i}})}\; .
$$
\end{theorem}

Let $S^{\C}_{C}$ be the usual semigroup of the complex branch $C$, i.\,e.,
$S^{\C}_{C} = \{\nu (f) : f\in \OO_{\C^2,0}\}$. The set of multiplicities
$\{m_{\sigma_0}, \ldots, m_{\sigma_g}\}$ is
the minimal set of
generators of $S^{\C}_{C}$. Let $e_i :=\gcd \{\seq{m}{\sigma_0}{\sigma_i}\}$ for
$i=0,1,\ldots, g$ and let $N_i = e_{i-1}/e_i$ for $i=1,2,\ldots,g$. (One
has
$m_{\tau_i}=N_i m_{\sigma_i}$ for $i=1,2,\ldots, g$: see, e.\,g., \cite{Wall} or \cite{Casas}.)

For $\sigma\in \Gamma$, let $\pi_{\sigma}: (X_{\sigma},D_{\sigma})\to (\C^2,0)$ be
the minimal modification of $(\C^2,0)$ such that $E_{\sigma}\subset D_{\sigma}$. In
particular, $E_{\sigma}$ is the last exceptional component appearing in $X_{\sigma}$
and is produced by blowing-up at a point $p_{\sigma}$ of a previous component.
For $f\in \OO_{\C^2,0}$, we will denote by
$e_{\sigma}(f)$ the multiplicity of the strict transform by $\pi_{\sigma}$ of the
curve $\{f=0\}$
at the point $p_{\sigma}$. Notice that $e_\sigma(f)$ coincides with the intersection
multiplicity of the strict transform of $\{f=0\}$ and $E_{\sigma}$ in the surface
$X_{\sigma}$.

For $j=1,\ldots, s$, let $\ell_j:= [\K_j : \K_{j-1}] =
[G_{j-1}:G_j] =\# (G_{j-1}/{G_j})$.

\begin{lemma}\label{lemma-1}
Let $\omega\in \check \Gamma$ be such that $\rho_d < \omega \le  \rho_{d+1}$ and
$\tau_{p} < \omega \le  \tau_{{p}+1}$.
One has
$$
M_\omega = \sum_{i=0}^{{p}} k_i M_{\sigma_i} +
 \sum_{j=1}^d r_j (\ell_j-1)M_{\rho_j}
$$
for some non negative integers $k_i$ $(0\le i \le {p})$, $k_i<N_i$ for $i\geq 1$,
and positive integers
$r_j$ $ (1\le j \le d)$.
\end{lemma}

\begin{proof}

Let us assume that  $\tau_{q-1}<\rho_d \le  \tau_q$
(see Figure \ref{fig4-1} below). For
$i=0,\ldots, g$, let $M_i := M_{\sigma_i}$, $m_i := m_{\sigma_i}$ and
let $\varphi_i=0$ be the equation of a curvette at the component
$E_{\sigma_i}$.

Let $\omega$ be such that
$\rho_d < \omega \le \rho_{d+1}$. Then one has
\begin{align*}
M_\omega &= m_{\omega} + ([G_0:G_d]-[G_1:G_d])e_{\rho_1}(\varphi_{\omega}) m_{\rho_1}
+ ([G_1:G_d]-[G_2:G_d])e_{\rho_2}(\varphi_{\omega}) m_{\rho_2} + \\
& \qquad \qquad + \cdots +  ([G_{d-1}:G_d]-1)e_{\rho_d}(\varphi_{\omega}) m_{\rho_d}
=
\\
&= m_\omega + \sum_{j=1}^d ([G_{j-1}:G_d]-[G_j:G_d])e_{\rho_j}(\varphi_\omega)
m_{\rho_j} =
m_\omega + \sum_{j=1}^d L^d_j e_{\rho_j}(\varphi_\omega)m_{\rho_j}
\end{align*}
where $L^d_j = [G_{j-1}:G_d]-[G_j:G_d]$, $1\le j \le d$.
In the same way, for $i=q,\ldots,{p}$ one has
\begin{equation}\label{eq-0}
M_i = m_i + \sum_{j=1}^d L^d_j e_{\rho_j}(\varphi_i)m_{\rho_j}\, .
\end{equation}



\begin{figure}[h]
$$
\unitlength=0.70mm
\begin{picture}(120.00,40.00)(0,10)

\thinlines
\put(-15,30){\circle*{2}}
\put(-25,24){{${\bf 1}=\sigma_0$}}

\put(-15,30){\line(1,0){5}}
\put(-8,30){\circle*{0.5}}
\put(-5,30){\circle*{0.5}}
\put(-2,30){\circle*{0.5}}
\put(1,30){\circle*{0.5}}
\put(4,30){\circle*{0.5}}
\put(5,30){\line(1,0){40}}

\put(10,30){\line(0,-1){20}}
\put(10,30){\circle*{2}}
\put(10,10){\circle*{2}}
\put(11.5,7){{$\sigma_{q-1}$}}
\put(9,33){{$\tau_{q-1}$}}

\put(30,30){\circle*{2}}
\put(27,23){{$\rho_d$}}

\put(30,30){\line(1,1){10}}
\put(30,30){\line(1,2){8}}
\put(30,30){\line(0,1){10}}

\put(40,30){\line(0,-1){20}}
\put(40,30){\circle*{2}}
\put(40,10){\circle*{2}}
\put(41.5,7){{$\sigma_{q}$}}
\put(39,33){{$\tau_{q}$}}

\put(48,30){\circle*{0.5}}
\put(51,30){\circle*{0.5}}
\put(54,30){\circle*{0.5}}

\put(55,30){\line(1,0){40}}

\put(60,30){\line(0,-1){20}}
\put(60,30){\circle*{2}}
\put(60,10){\circle*{2}}
\put(61.5,7){{$\sigma_{{p}}$}}
\put(59,33){{$\tau_{{p}}$}}

\put(75,30){\circle*{3}}
\put(74,23){{$\omega$}}

\put(90,30){\line(0,-1){20}}
\put(90,30){\circle*{2}}
\put(90,10){\circle*{2}}
\put(91.5,7){{$\sigma_{{p}+1}$}}
\put(89,33){{$\tau_{{p}+1}$}}

\put(98,30){\circle*{0.5}}
\put(101,30){\circle*{0.5}}
\put(104,30){\circle*{0.5}}
\put(107,30){\circle*{0.5}}

\put(112,30){\line(1,0){20}}
\put(120,30){\circle*{2}}
\put(118,23){{$\rho_{d+1}$}}

\put(120,30){\line(2,1){10}}
\put(120,30){\line(1,1){10}}
\put(120,30){\line(0,1){10}}

\put(134,30){\circle*{0.5}}
\put(137,30){\circle*{0.5}}
\put(140,30){\circle*{0.5}}

\end{picture}
$$
\caption{The position of the vertices.}
\label{fig4-1}
\end{figure}

Moreover, it is known  that
$m_\omega\in \langle\seq{m}0{p}\rangle$, where
$\langle\seq{m}0{p}\rangle$ stands for the subsemigroup of $\Z_{\ge 0}$ generated by
$\seq{m}0{p}$, i.e.
$$
\langle \seq m0p \rangle := \{n\in \N : n = a_0 m_0 + \cdots + a_p m_p \text{ with }
a_i\in \Z_{\ge 0}
\}\; ,
$$
and so
$m_\omega = k_0 m_0 + \sum_{i=1}^{p} k_i m_i$, for some $k_0\ge 0$ and
$0\le k_i <N_i$ for $i\ge 1$ (which are unique with these restrictions).
Then one can write
\begin{equation}\label{eq-1}
M_\omega =
\sum_{0}^{p} k_i m_i + \sum_{j=1}^d L^d_j e_{\rho_j}(\varphi_\omega)m_{\rho_j}
\end{equation}

Let us assume that $q<{p}$, note that in this case we have
$\tau_{q-1}<\rho_d \le \tau_{q}<\tau_p \le \omega$ (see Figure \ref{fig4-1}).
For $j\le d$
one has
$e_{\rho_j}(\varphi_{\omega})/e_{\tau_{{p}-1}}(\varphi_\omega) =
e_{\rho_j}(\varphi_{{p}})/e_{\tau_{{p}-1}}(\varphi_{p})$, moreover
$e_{\tau_{{p}-1}}(\varphi_{p})=1$, so
$e_{\rho_j}(\varphi_{\omega})=
e_{\tau_{{p}-1}}(\varphi_{\omega}) e_{\rho_j}(\varphi_{p})$.
It is known that $e_i = e_{\tau_i}(\varphi_C)$, $i\le g$, where $\varphi_C=0$ is an
equation of the complex curve $C$. As a consequence, for $i\le p$, one has
$N_i = e_{i-1}/e_i = e_{\tau_{i-1}}(\varphi_\omega)/e_{\tau_i}(\varphi_\omega)$ and
so
\begin{equation}\label{eq-1-2}
e_{\rho_j}(\varphi_\omega) = N_{p} e_{\tau_{p}}(\varphi_\omega)e_{\rho_j}(\varphi_{p})\; .
\end{equation}
For $i,j$ such that $\tau_i > \tau_{i-1}>\rho_j$, the
same computations applied to
$\varphi_i$ (instead of $\varphi_\omega$) and
$\varphi_{i-1}$ (instead of $\varphi_p$)
show that
\begin{equation}\label{eq-1-3}
e_{\rho_j}(\varphi_i) = N_{i-1} e_{\tau_{i-1}}(\varphi_i) e_{\rho_j}(\varphi_{i-1})
= N_{i-1} e_{\rho_j}(\varphi_{i-1})\; .
\end{equation}

Using equality (\ref{eq-1-2}) in the right part of equation (\ref{eq-1}),
one gets:

\begin{align*}
\sum_{j=1}^d L^d_j e_{\rho_j}(\varphi_\omega)m_{\rho_j} &= N_{p}
e_{\tau_{p}}(\varphi_\omega)
\sum_{j=1}^d L^d_j e_{\rho_j}(\varphi_{p})m_{\rho_j} = \\
& = k_{p} \sum_{j=1}^d L^d_j e_{\rho_j}(\varphi_{p})m_{\rho_j} +
k'_{p} \sum_{j=1}^d L^d_j e_{\rho_j}(\varphi_{p})m_{\rho_j}
\end{align*}
with $k'_{p} = N_{p} e_{\tau_{p}}(\varphi_\omega)-k_{p}\ge
N_{p} - k_{p} \ge 1$,
with $k_p$ being the integer appearing in equation (\ref{eq-1}).
Now, using the expression of $M_i$ in (\ref{eq-0}) for $i=p$, one has
$$
M_\omega =
\sum_{0}^{{p}-1}k_i m_i +
k_{p} M_{p} +
k'_{p} \sum_{j=1}^d L^d_j e_{\rho_j}(\varphi_{p})m_{\rho_j}\; .
$$
In the same way, using repeatedly that
$e_{\rho_j}(\varphi_i)=N_{i-1}e_{\rho_j}(\varphi_{i-1})$
and the equation (\ref{eq-0}) for $M_i$, $i\ge q$,
one gets:
\begin{equation}\label{eq-2}
M_\omega =
\sum_{i=0}^{q-1}k_i m_i +
\sum_{i=q}^{p} k_i M_i +
r_d \sum_{j=1}^d L^d_j e_{\rho_j}(\varphi_q)m_{\rho_j}
\end{equation}
for some integer $r_d\ge 1$.
Now, one can use
the equality $L^d_{j} = [G_{d-1}:G_d]
L^{d-1}_j$ for
$j\le d-1$ and (\ref{eq-1-3})
in the last summand of (\ref{eq-2}) to get:
\begin{equation}\label{eq-2-2}
r_d
\sum_{j=1}^d L^d_j e_{\rho_j}(\varphi_q)m_{\rho_j} =
 r_d
[G_{d-1}:G_d]\sum_{j=1}^{d-1} L^{d-1}_j e_{\rho_j}(\varphi_q)m_{\rho_j} +
r_d
L^d_d e_{\rho_d}(\varphi_q)m_{\rho_d}\; .
\end{equation}

Using that $L^d_d = [G_{d-1}:G_d]-1$,
$e_{\rho_j}(\varphi_{q})=e_{\rho_j}(\varphi_{\rho_d})$,
$e_{\rho_d}(\varphi_{\rho_d})=1$ and the corresponding expression for $M_{\rho_d}$ in
the same terms as we use at the beginning of the proof for $M_\omega$, the
equality (\ref{eq-2-2}) can be written as:
\begin{align*}
& r_d ([G_{d-1}:G_d]-1)\left(\sum_{j=1}^{d-1} L^{d-1}_j e_{\rho_j}(\varphi_{\rho_d})
m_{\rho_j} + m_{\rho_d}\right) +
r_d \sum_{j=1}^{d-1}L^{d-1}_{j}e_{\rho_j}(\varphi_q)m_{\rho_j} = \\
& =
r_d ([G_{d-1}:G_d]-1)M_{\rho_d} +
r_d \sum_{j=1}^{d-1}L^{d-1}_{j}e_{\rho_j}(\varphi_q)m_{\rho_j}\; .
\end{align*}
Thus, using this equality and taking into account that
$\ell_j = [G_{j-1}:G_j]$, the equation
(\ref{eq-2}) for $M_{\omega}$ becomes:

\begin{equation}\label{eq-3}
M_\omega =
\sum_{i=q}^{p} k_i M_i +
r_d (\ell_d-1)M_{\rho_d} +
\left(\sum_{0}^{q-1}k_i m_i  +
r_d \sum_{j=1}^{d-1}L^{d-1}_{j}e_{\rho_j}(\varphi_q)m_{\rho_j}
\right)\; .
\end{equation}
If $d=1$, one has $m_i=M_i$ for $i=0,\ldots, q-1$ and the proof is finished.
Otherwise,
we repeat the above computations
for the last two summands in (\ref{eq-3})
(in brackets):
\begin{equation}\label{eq3-1}
\sum_{0}^{q-1}k_i m_i  +
r_d \sum_{j=1}^{d-1}L^{d-1}_{j}e_{\rho_j}(\varphi_q)m_{\rho_j}
\; .
\end{equation}
To do that, note that
the equation (\ref{eq-0}) remains valid for $M_i$, $i\le q$, and that
the equation (\ref{eq3-1}) is essentially the same as the one in the
right hand side of (\ref{eq-1}):
\begin{enumerate}
\item The factor
$r_d\ge 1$ does not affect the computations in the proof.
\item Equation (\ref{eq-1-3}) could be applied to $\varphi_{p}$ (and to the
succesive $\varphi_i$, $i<p$) instead of equation (\ref{eq-1-2}) for
$\varphi_{\omega}$.
\item The succesive integers $k'_i \ge N_i -k_i$ are always greater or equal to 1.
\end{enumerate}


As a consequence, we get:
\begin{align*}
M_\omega
& = \sum_{i: \rho_{1}<\tau_i\le\tau_{p}} k_i M_i +
\sum_{i: \tau_i\le \rho_{1}} k_i m_i
+ \sum_{j=1}^d r_j (\ell_j-1)M_{\rho_j} \\
& = \sum_{i=0}^{{p}} k_i M_i +
 \sum_{j=1}^d r_j (\ell_j-1)M_{\rho_j}
\end{align*}
for some integers $r_j\ge 1$ ($1\le j\le d$).
\end{proof}

\begin{proposition}\label{prop:semig}
The set of integers $\{\seq{M}{\sigma_0}{\sigma_g}\}$ is a set of
generators of  the semigroup of values
$S^{\K}_C$ of the valuation $\nu$ on $\EE_{\K^2,0}$.
\end{proposition}


\begin{proof}
It is clear that the semigroup $S^{\K}_{\nu}$ is generated by the integers
$M_{\sigma}$, for all $\sigma\in \check \Gamma$. Thus one has to show that
$M_\sigma\in \langle M_{\sigma_0}, \ldots, M_{\sigma_g}\rangle$ for any $\sigma$.

Let $[\alpha,\beta]$ denote the geodesic
in $\check \Gamma$ between the vertices $\alpha$ and $\beta$.
If $\omega\in [\sigma_0, \delta_C]$, Lemma \ref{lemma-1} gives
$$
M_\omega = \sum_{i=0}^{h} k_i M_{\sigma_i} +
 \sum_{j=1}^d r_j (\ell_j-1)M_{\rho_j} \; ,
$$
where $\rho_d < \omega \le  \rho_{d+1}$, $\tau_h < \omega \le  \tau_{h+1}$,
$k_i$ are non-negative integers and $r_j$ are positive integers.
The same equation applied to $\rho_j$ gives us that
$M_{\rho_j}\in \langle \seq M0{i}\rangle$ if
$\tau_i \le \rho_j < \tau_{i+1}$. In particular,
$M_{\rho_j}\in \langle \seq M0{h}\rangle \subset
\langle \seq M0{g}\rangle$

\medskip

If $\omega\in [\sigma_i, \tau_{i}]$ for some $i =1,\ldots, g$, then
$\nu(\varphi_\omega) = k m_{\sigma_i}$ for a positive integer $k$.
For the
conjugate $g \omega$, $g\in G$,
one has that
$g{\omega}\in [g{\sigma}_i, g{\tau}_{i}]$ and so
$\nu(g {\varphi_\omega}) = \nu_{g {C}}(\varphi_\omega)= k m_{g {\sigma}_i}$. Thus, in
this case one has that
$M_\omega = \nu(\sum g \varphi_\omega) = k(\sum m_{g \omega}) =
k M_{\sigma_i}$. Note that
$M_{\tau_i} = N_i M_{\sigma_i}$.

\end{proof}

\begin{remark}\label{rmk5}
For $\sigma\in \check \Gamma$, $\tau_k <  \sigma \le  \tau_{k+1}$, it is known  that
$m_\sigma\in \langle m_{\sigma_0}, \ldots, m_{\sigma_k} \rangle$ or equivalently
that $e_k$ divides $m_{\sigma}$ (see \cite{Wall}).
Taken into account the expression for $M_i$
given in equation (\ref{eq-0}),
this implies that
$e_{i-1}$ divides
$\sum_{j=1}^d L^d_j e_{\rho_j}(\varphi_i)m_{\rho_i}$ and so
$\gcd (\seq M0i) = \gcd (\seq m0i) = e_i$ and
$$N_i = e_{i-1}/e_i = \gcd(\seq M0{i-1})/\gcd(\seq M0i)\, . $$
\end{remark}

\begin{proposition}\label{lemma-2}
For $i=1,\ldots , g$,  one has
 $(N_{i}-1)M_{\sigma_i}\notin \langle M_{\sigma_0}, \ldots, M_{\sigma_{i-1}}\rangle$,
 $N_{i} M_{\sigma_{i}}\in \langle M_{\sigma_0}, \ldots, M_{\sigma_{i-1}}\rangle$
and
$N_{i-1} M_{\sigma_{i-1}}< M_{\sigma_{i}}$.
In particular, $\{M_{\sigma_0}, \ldots, M_{\sigma_g}\}$ is the minimal set of
generators of the semigroup $S^{\K}_C := \langle M_{\sigma_0}, \ldots, M_{\sigma_g}
\rangle$ and,
moreover, $S^{\K}_C$ is the semigroup of values of a complex analytic irreducible
curve.
\end{proposition}

\begin{proof}
If $(N_i-1)M_i\in \langle \seq M0{i-1}\rangle$ then,
as a consequence of Remark \ref{rmk5}, $e_{i-1}$ divides
$(N_i-1)M_i$ and so $e_{i-1}$ divides
$(N_i-1)m_i$. But this contradicts the properties of the minimal set of
generators of the semigroup $S^{\C}_{C}$.

Note that $N_i M_i = M_{\tau_i}$, then
$N_i M_i \in \langle \seq M0{i-1}\rangle$ as a consequence of Proposition
\ref{prop:semig} (more precisely of the first part of its proof).

Let us prove that, for $i=1,\ldots, g$, one has that $N_{i-1} M_{i-1}< M_i$.
If $\rho_{d} \le \tau_{i-1} < \tau_i$, the equation (\ref{eq-0}) for $M_i$ together with
the equality $e_{\rho_i} = N_{i-1}e_{\rho_i}(\varphi_{i-1})$ for $i\le d$ imply
that
\begin{align*}
M_i & = m_i + \sum_{j=1}^d L^d_j e_{\rho_j}(\varphi_i)m_{\rho_j} =
m_i + N_{i-1}\sum_{j=1}^d L^d_j e_{\rho_j}(\varphi_{i-1})m_{\rho_j}  > \\
& > N_{i-1}\left( m_{i-1} + \sum_{i=j}^d L^d_j
e_{\rho_j}(\varphi_{i-1})m_{\rho_j}\right)
= N_{i-1}M_{i-1}\; .
\end{align*}

If $\tau_{i-1}< \rho_i \le  \tau_i$ the equation (\ref{eq-0}) for $M_{i-1}$ is
$M_{i-1}  = m_{i-1} + \sum_{j=1}^{d'} L^{d'}_j e_{\rho_j}(\varphi_{i-1})m_{\rho_j} $
for some $d'<d$. Thus, using the equality
$$
L^d_{j} = [G_{d-1}:G_d] L^{d-1}_j = \ell_{d-1}L^{d-1}_j = \ell_{d-1}\cdots
\ell_{d'}L^{d'}_j \; ,
$$
one gets
\begin{align*}
M_i & = m_i + \sum_{j=1}^d L^d_j e_{\rho_j}(\varphi_i)m_{\rho_j}
> N_{i-1}\left( m_{i-1} + \sum_{j=1}^d L^d_j
e_{\rho_j}(\varphi_{i-1})m_{\rho_j}\right)
\\
& > N_{i-1}\left( m_{i-1} + \sum_{j=1}^{d'} L^{d'}_j
e_{\rho_j}(\varphi_{i-1})m_{\rho_j}\right) =
N_{i-1}M_{i-1}\; .
\end{align*}

It is known that the proved properties implies that the semigroup
$S = \langle M_{\sigma_0}, \ldots, M_{\sigma_g}  \rangle$ is the semigroup of a
complex irreducible plane curve (see \cite{Bresinsky}).

\end{proof}

\begin{Proof} {\bf of Theorem 1.}
The statement follows directly from
Propositions~\ref{prop:semig} and
\ref{lemma-2} taking into
account that
$M_{\tau_i}=N_i M_{\sigma_i}$ for all $i\ge 1$.
\end{Proof}

\begin{remark}
As we have proved, $S^\K_{\nu}$ is the semigroup of a complex analytic branch $C'$ (and so
$P^S_C(t)$ is the usual Poincar\'e series of it) different from $C$ if
$m_{\rho_1}<m_{\tau_g}$.
We do not know a geometric reason for the existence of the curve $C'$.
\end{remark}

\section{The classical Poincar\'e series of a curve valuation}\label{sect:PS-curve}
\begin{theorem}\label{theo:theo2}
 In Cases I and II, one has
 \begin{eqnarray*}
  P^{\K}_{C}(t)&=&P_{C}^S(t)\cdot
\prod_{j=1}^s \frac{1-t^{\ell_j M_{\rho_j}}}{1-t^{M_{\rho_j}}}
\\
  {\ }&=&\frac{
\prod_{i=1}^g(1-t^{M_{\tau_i}})}{\prod_{i=0}^g(1-t^{M_{\sigma_i}})}\cdot
\prod_{j=1}^s \frac{1-t^{\ell_j M_{\rho_j}}}{1-t^{M_{\rho_j}}}
 \end{eqnarray*}
 ($s=\infty$ in Case II).
\end{theorem}

Before proving the Theorem, we elaborate some statements.

 Let $P^{\K}_C(t)=\sum_{v=0}^{\infty} a_vt^v$. For a collection
 $\{k_{\sigma}\}$,
 $\sigma\in\check{\Gamma}$, $k_{\sigma}\in\Z_{\ge0}$, let
 $\calE_{\K^2,0}^{\{k_{\sigma}\}}$ be the set of germs $f\in\calE_{\K^2,0}$
 such that the intersection multiplicity of the strict transform of the zero-level
set
 $\{f=0\}$ with the component $E_{\sigma}$ of the exceptional divisor $D$ is equal
 to $k_{\sigma}$ and, moreover, this strict transform intersects $D$ only at
 smooth points of the total transform $\pi^{-1}(C)$ of the curve $C$.
 (Pay attention that the intersection multiplicity of the strict transform of
 $\{f=0\}$ with the component $E_{g\sigma}$,
 $g\in G$, is also equal
 to $k_{\sigma}$.)
 Let $\nu(\{k_{\sigma}\}):=\sum_{\sigma\in\check{\Gamma}}k_{\sigma}M_{\sigma}$. One
can see that $\nu(\{k_{\sigma}\})=\nu(f)$
 for any $f\in\calE_{\K^2,0}^{\{k_{\sigma}\}}$.

 \begin{remark}
 Without loss of generality, we shall use the assumption that, for a fixed
 $v=\nu(\{k_{\sigma}\})$, the condition
that the zero-level set $\{f=0\}$ with $f\in \calE_{\K^2,0}$ intersects $D$ only at
smooth points holds automatically. This can be achieved if one makes sufficiently
many additional
blow-ups at intersection points of the components of the total transform
$\pi^{-1}(C)$.
 \end{remark}

 Let
 $F^{\{k_{\sigma}\}}$ be the image of
 $\calE_{\K^2,0}^{\{k_{\sigma}\}}$ in
 the quotient $J^{\C}(v)/J^{\C}(v+1)\cong\C$ with
 $v=\nu(\{k_{\sigma}\})$.(Of course
 $F^{\{k_\sigma\}}\subset J^{\K}(v)/J^{\K}(v+1) \subset
 J^{\C}(v)/J^{\C}(v+1)$.)
 Let $\overline{F}^{\{k_{\sigma}\}}$ be the linear span
 over $\K$ of $F^{\{k_{\sigma}\}}$.

 One has:
 \begin{enumerate}
  \item[1)] $$\bigcup_{\{k_{\sigma}\}: \nu(\{k_{\sigma}\})=v}
  F^{\{k_{\sigma}\}}=
  J^{\K}(v)/J^{\K}(v+1) \setminus \{0\};$$
  \item[2)] for each collection $\{k_{\sigma}\}$, $F^{\{k_{\sigma}\}}$
  is the complement to a finite hyperplane arrangement
  in the vector $\K$-space $\overline{F}^{\{k_{\sigma}\}}$.
 \end{enumerate}

Let $d^{\{k_{\sigma}\}}$ be the dimension (over $\K$) of the vector space
$\overline{F}^{\{k_{\sigma}\}}$.
One obviously has
$a_v=0$ if there are no $\{k_{\sigma}\}$ with $\nu(\{k_{\sigma}\})=v$,
$$
a_v=\max_{\{k_{\sigma}\}:\nu(\{k_{\sigma}\})=v}d^{\{k_{\sigma}\}}
{\text{\ \ \ otherwise}}\,.
$$

For $j=1,2, \ldots, s+1$ ($s$ may be equal to $+\infty$), let $A_j$ be
the set of collections $\{k_{\sigma}\}$ such that $k_{\sigma}=0$
for all $\sigma$
on the geodesic $[\rho_j,\delta_C]$ from $\rho_j$ to $\delta_C$ (including the ends).
For $j=s+1$ (if $s<+\infty$), we assume this geodesic to be empty.
Let the series
$P^{(j)}_C(t)=\sum_{v=0}^{\infty} a^{(j)}_v t^v$ be defined by
$$
a^{(j)}_v=\max_{\{k_{\delta}\}\in A_j:\nu(\{k_{\sigma}\})=v}d^{\{k_{\sigma}\}}.
$$
One obviously has $P^{\K}_C(t)=P^{(s+1)}_C(t)$ if $s<+\infty$ and, if $s=+\infty$
(i.\,e., in Case II : the process of the $G$-resolution of the curve $C$ is
infinite),
$P^{\K}_C(t)$ is the limit of $P^{(j)}_C(t)$
for $j\to+\infty$ in ``the $\mathfrak{m}$-adic topology''
($\mathfrak{m}$ is the maximal ideal in $\Z[t]$).
The latter means that, for any $v$, $a_v=a_v^{(j)}$ for $j$ large enough.

For a function $f\in\calE_{\K^2,0}$,
let $\widetilde{f}=f\circ\pi$ be the lifting of $f$ to the surface of the resolution.
In order to compute (or evaluate) the dimension $d^{\{k_{\sigma}\}}$, we shall
regularly use the following construction exploited first in~\cite{IJM}.
Let $h\in \calE_{\K^2,0}$ be a {\bf fixed
function} with $\nu(h)=\nu(\{k_{\sigma}\})$.
(In fact we shall take a fixed function from
$\calE_{\K^2,0}^{\{k_{\sigma}\}}$.)
For each function
$f\in\calE_{\K^2,0}^{\{k_{\sigma}\}}$
we shall consider the ratio
$\Psi={\widetilde f}/{\widetilde h}$.
Let $q$ of the exceptional divisor with the strict transform
of the curve $C$.
The function $\Psi$ is regular in a neighbourhood of $q$ and its value
$\Psi(q)$ is different from zero.
The set $\left\{({\widetilde f}/{\widetilde h})(q):
f\in\calE_{\K^2,0}^{\{k_{\sigma}\}}\right\}$ of values of these ratios at the point
$q$ as a subset of $\C$ is obtained from $F^{\{k_{\sigma}\}}$, also considered
as a subset of $\C\cong J^{\C}(v)/J^{\C}(v+1)$ with
$v=\nu(\{k_{\sigma}\})$, by multiplication by a non-zero
factor.
This follows from the following arguments.
The set $F^{\{k_{\sigma}\}}$ is just the set of
the coefficients $a(f)$ in Equation~(\ref{eqn:curve_valuation}) for
$f\in \calE_{\K^2,0}^{\{k_{\sigma}\}}$.
The coefficient $a(f)$ is equal to
the coefficient $a(h)$ (a fixed number)
multiplied by the value of the ratio $\Psi$.

\medskip

Assume first that $\rho_1\ne\sigma_0$.
The necessary changes in the case $\rho_1=\sigma_0$
will be discussed later.

 \begin{lemma}\label{lemma:2_in_Th2}
  $a^{(1)}_v\ne 0$ if and only if $v\in S^{\K}_C$.
 \end{lemma}

 \begin{proof}
 Obviously $a^{(1)}_v=0$ for
 $v\not\in S^{\K}_C$.
  Any $v\in S^{\K}_C$ can be represented
  in the form $v=\sum_{i=0}^g k_iM_{\sigma_i}$
  (Proposition~\ref{prop:semig}).
  Any function $f\in\EE_{\K^2,0}^{\{k_{\sigma}\}}$ with
  $k_{\sigma_i}=k_i$, $k_{\sigma}=0$ otherwise
  gives a non-zero element in the corresponding
 $F^{\{k_{\sigma}\}}$. Therefore $d^{\{k_{\sigma}\}}\ge 1$
 and $a^{(1)}_v\ge 1$.
 \end{proof}

 \begin{lemma}\label{lemma:1_in_Th2}
  For each $v$, $a^{(1)}_v\le 1=[\K_0:\K_0]$.
  As a consequence, $a^{(1)}_v= 1=[\K_0:\K_0]$
  if and only if $v\in S^{\K}_C$.
 \end{lemma}

 \begin{proof}
 Let $\{k_{\sigma}\}\in A_1$ and let $h$ be a fixed function from
$\calE_{\K^2,0}^{\{k_{\sigma}\}}$.
 For $f\in\calE_{\K^2,0}^{\{k_{\sigma}\}}$, the ratio
 $\Psi={\widetilde f}/{\widetilde h}$ has no zeroes or poles
 on all the components from the union of the geodesics
$[\rho_1, g \delta_C]$ for all $g\in G=G_0$. Therefore
 this is a (non-zero) constant on them. This constant is invariant with respect to
 $G$ and therefore belongs to $\K$.
 \end{proof}

Lemmas~\ref{lemma:2_in_Th2} and~\ref{lemma:1_in_Th2} obviously imply the following
statement.

 \begin{proposition}\label{prop:Prop1_in_Th2}
 One has
 $$
 P^{(1)}_C(t)=P^{S}_C(t)=\frac{\prod_{i=1}^g (1-t^{M_{\tau_i}})}
{\prod_{i=0}^g (1-t^{M_{\sigma_i}})}\;.
 $$
 \end{proposition}

 \begin{lemma}\label{lemma:3_in_Th2}
  For each $j\ge 1$ and any $v$, $a^{(j)}_v\le [\K_{j-1}:\K_0]$.
 \end{lemma}

 \begin{proof}
 The proof is essentially the same as the one of Lemma~\ref{lemma:1_in_Th2}. For
 $\{k_{\sigma}\}\in A_j$ and
 a fixed function
 $h\in\calE_{\K^2,0}^{\{k_{\sigma}\}}$,
 the ratio
 $\Psi={\widetilde f}/{\widetilde h}$
 ($f\in\calE_{\K^2,0}^{\{k_{\sigma}\}}$)
 has no zeroes or poles
 on all the components from the union of the geodesics
$[\rho_j, g\delta_C]$ for all $g\in G_{j-1}$. Therefore
 this is a (non-zero) constant on them. This constant is invariant with respect to
$G_{j-1}$ and therefore belongs to $\K_{j-1}$.
 \end{proof}

 \begin{lemma}\label{lemma:4_in_Th2}
  Let $v\in S^{\K}_C$ be of the form
  $v=\ell M_{\rho_j}+b$ with $b\in S^{\K}_C$
  and such that $b-M_{\rho_j}\notin S_C^{\K}$.
  Assume that $\ell\le \ell_j-1$ $(\ell_j=[\K_j:\K_{j-1}])$.
  We also suppose that $v$ is not of the form $M_{\varkappa}+b'$ with
  $\varkappa>\rho_j$ and $b'\in S^{\K}_C$.
  (The corresponding case will be analyzed in
  Lemma~\ref{lemma:6_in_Th2}.)
  Then
  $a_{v}^{(j+1)}=
  \ell[\K_{j-1}:\K_0]+a_{b}^{(j)}$.
 \end{lemma}

\begin{proof}
 We shall use the coordinate $w$ on
 $E_{\rho_j}$
 (from the corresponding standard coordinate chart) such that $w=\infty$ at the intersection
 point with the previous component on the
 the geodesic between $\sigma_0$ and $\rho_j$.
 Let $w=z_0\in \K_j\setminus\K_{j-1}$ be
 a point on $E_{\rho_j}$ blown-up in the process of the resolution of the curve $C$.
 One has $\K_j = \K_{j-1}[z_0]$
 and
$1$, $z_0$, $z_0^2$, \dots, $z_0^{\ell_j-1}$
is a basis of the vector space $\K_j$
over $\K_{j-1}$.

 Let us represent a collection $\{k_{\sigma}\}$
 with $\nu(\{k_{\sigma}\})=v$ in the form
 $\{k_{\sigma}\}=\{k^1_{\sigma}\}+\{k^2_{\sigma}\}$,
 where $k^1_{\rho_j}=\ell'$ ($\ell'\le\ell$),
 $k^1_\sigma=0$ for $\sigma\ne \rho_j$,
 $k^2_{\rho_j}=0$.
 A function $f\in\calE_{\K^2,0}^{\{k_{\sigma}\}}$
 can be represented as $f=f_1f_2$ with
 $f_r\in\calE_{\K^2,0}^{\{k^r_{\sigma}\}}$, $r=1,2$.
 Let $h_r$, $r=1,2$, be fixed functions from
 $\calE_{\K^2,0}^{\{k^r_{\sigma}\}}$.
 For certainty, as $h_1$ we shall take
 the $\ell'$-th power of the left hand side
 of the equation of a $G_{j-1}$-curvette at the point $w=1$ of the component $E_{\rho_j}$.
 Let $h=h_1h_2$.
 One has $\Psi={\widetilde f}/{\widetilde h}=
 \Psi_1\Psi_2$, where
 $\Psi_r={\widetilde f_r}/{\widetilde h_r}$.
  The function $\Psi_2$ is constant on the union of all the components $E_\sigma$
 for $\sigma$ from the
 geodesics between $\rho_j$ and
 $g\delta_C$ with $g\in G_{j-1}$ (including the both end points of the geodesic).
 Its value on them belongs to $\K_{j-1}$
 and the maximal possible dimension (over $\K$)
 of the set of its values is $a_{b}^{(j)}$.
 The function $\Psi_1$ is constant on the union of all the components $E_\sigma$
 for $\sigma$ from the
 geodesic between $\rho_j$ and
 $g \delta_C$ with $g\in G_j$ {\bf excluding}
 $\rho_j$ itself.
 Its restriction to the component $E_{\rho_j}$
 is equal to
 $$
\Psi_1 |_{E_{\rho_j}} =
\frac{ c\left(\prod\limits_{p=1}^{\ell'}(w-w_p)\right)}
{(w-1)^{\ell'}}\, ,
 $$
 where $w_1,\ldots, w_{\ell'}$ are the intersection points
 of the strict transform of the curve
 $\{f_1=0\}$
 with the component $E_{\rho_j}$.
 All these intersection points $w_p$
 belong to $\K_{j-1}$. This follows
 from the fact that the set of intersection points of the strict transform of the curve
 $\{f_1=0\}$ with $E_{\rho_j}$ is
 $G_j$-invariant, whereas the orbit of a point
 not belonging to $\K_{j-1}$ consists of
 $\ell_j$ points.
 Moreover, its value at $w=\infty$,
 that is
 $c$, belongs to $\K$. (This follows from the fact that this ratio is constant on the union of
 the components from the geodesic
 $[\sigma_0, \rho_j)$.)
 Thus up to the constant $(z_0-1)^{\ell'}$
 (which does not influence the dimensions) the value of $\Psi_1$ on the union of the
components from the geodesic
 $(\rho_j, \delta_C]$ is constant and is
 equal to
$c\prod_{p=1}^{\ell'}(z_0-w_p)=
 c(z_0^{\ell'}-(w_1+\ldots+w_{\ell'})z_0^{\ell'-1}+\ldots
 \pm (w_1\cdot\ldots\cdot w_{\ell'}))$.
 Therefore the value of the ratio $\Psi$ on this union
 is equal to $c(\alpha_{\ell'}z_0^{\ell'}+
 \alpha_{\ell'-1}z_0^{\ell'-1}+\alpha_{1}z_0+\alpha_0)$,
 where $\alpha_0$, \dots, $\alpha_{\ell'-1}$ are generic points of $\K_{j-1}$
(because we can choose arbitrarily $\seq w1{\ell'}$ in $\K_{j-1}$)
and
$\alpha_{\ell'}$ (also a point of
 $\K_{j-1}$) is the value
 of the ratio $\Psi_2$. The (maximal) dimension of the set of values of $\alpha_{\ell'}$
is $a_{b^*}^{(j)}$, where
$b^*=v-\ell'M_{\rho_j}$.
 Since $1$, $z_0$, \dots $z_0^{\ell'}$ are linear independent over $\K_{j-1}$, the
total dimension is equal to $\ell'[\K_{j-1}:\K_0]+a_{b^*}^{(j)}$.
Thus it is maximal only for $\ell'=\ell$
 and is equal to
 $\ell[\K_{j-1}:\K_0]+a^{(j)}_{b}$.

 \end{proof}

 \begin{lemma}\label{lemma:5_in_Th2}
  If $\rho_j<\varkappa\le \delta_C$ (that is $\varkappa$ is
  on the geodesic between $\rho_j$ and $\delta_C$ and $\varkappa \neq \rho_j$; in particular
  for $\varkappa= \rho_{j+1}$), then
  $a_{M_{\varkappa}}^{(j+1)}=[\K_{j}:\K_0]$.
 \end{lemma}

 \begin{proof}
 Any collection $\{k_{\sigma}\}$ with $\nu(\{k_{\sigma}\})=M_{\varkappa}$
 obviously belongs to $A_{j+1}$.
 According to Lemma~\ref{lemma:3_in_Th2},
 $a_{M_{\varkappa}}^{(j+1)}\le[\K_{j}:\K_0]$.
 Let us take the collection $\{k^0_{\sigma}\}$ with
 $k^0_{\varkappa}=1$ and $k^0_{\sigma}=0$
 for $\sigma\ne \varkappa$. One has
 $\nu(\{k^0_{\sigma}\})=M_{\varkappa}$.
 Let $f=0$ be an equation of a $G_\varkappa(=G_j)$-curvette intersecting
 the component $E_{\varkappa}$ at the point with the coordinate
 $\alpha\in \K_j$. Here we assume that the coordinate is equal to
 zero at the intersection point with the component
 on the geodesic to $\delta_C$ and to $\infty$ at the one with
 the component on the geodesic to $\sigma_0$.
 Let $g$ be one of these functions, say, with $\alpha=1$. The ratio
 $\Psi={\widetilde f}/{\widetilde g}$, for $f\in \EE^{\{k^0_\sigma\}}_{\K^2,0}$,
has no zeroes or poles
 on all the components of each of the indicated geodesics and thus
 is constant on each of them. On the component $E_\varkappa$ this ratio is
 of the form $c\frac{w-\alpha}{w-1}$ with $c\in \K$
 and (arbitrary!) $\alpha\in \K_j$. The value on $\delta_C$ (equal to the value at
$w=0$, i.\,e., $c \alpha$) is an arbitrary point of $\K_{j}$.
 Therefore $d^{\{k^0_{\sigma}\}} =[\K_{j}:\K_0]$
 and $a_{M_{\varkappa}}^{(j+1)}\ge[\K_{j}:\K_0]$.
 \end{proof}

 \begin{lemma}\label{lemma:6_in_Th2}
  Let $\varkappa$ be as in Lemma~\ref{lemma:5_in_Th2}.
  Due to Lemma~\ref{lemma-1}, one has that
  $M_{\varkappa}=(\ell_j-1)M_{\rho_j}+b$
  with $b\in S^{\K}_C$. Then
  $a_{b}^{(j)}=[\K_{j-1}:\K_0]$.
 \end{lemma}

 \begin{proof}
 Due to Lemma~\ref{lemma-1}, one has
 $M_{\varkappa}=b'+\sum_{q=1}^j
 (\ell_q-1)M_{\rho_q}$ with $b'\in S^{\K}_C$.
 Let $c_k:=b'+\sum_{q=1}^k(\ell_q-1)M_{\rho_q}$,
 $k=1,\ldots, j$.
 We have to show that $a^{(j)}_{c_{j-1}}=
 [\K_{j-1}:\K_0]$.
 We shall show that $a^{(k)}_{c_{k-1}}=
 [\K_{k-1}:\K_0]$ for all $k$. We shall use induction in $k$. For $k=1$ this is
trivial:
 $a^{(1)}_{c_0}=1=[\K_{0}:\K_0]$.
 Due to Lemma~\ref{lemma:4_in_Th2},
 we have $a^{(k)}_{c_{k-1}}=(\ell_{k-1}-1)
 [\K_{k-2}:\K_0]+a^{(k-1)}_{c_{k-2}}$.
 By the induction assumption
 $a^{(k-1)}_{c_{k-2}}=[\K_{k-2}:\K_0]$
 and therefore
 $a^{(k)}_{c_{k-1}}=\ell_{k-1}[\K_{k-2}:\K_0]=
 [\K_{k-1}:\K_0]$.
 \end{proof}

 Lemmas~\ref{lemma:3_in_Th2}--\ref{lemma:6_in_Th2}
 imply the following statement.

 \begin{proposition}\label{prop:Prop2_in_Th2}
 One has
  \begin{eqnarray}
  P^{(j+1)}(t)&=&P^{(j)}(t)\cdot(1+t^{M_{\rho_j}}+
  t^{2M_{\rho_j}}+\ldots+t^{(\ell_j-1)M_{\rho_j}})
  \label{eqn:Prop2}
  \\
  &=&P^{(j)}(t)\cdot\frac{1-t^{\ell_j M_{\rho_j}}}{1-t^{M_{\rho_j}}}.\nonumber
  \end{eqnarray}
 \end{proposition}

 \begin{proof}
 Equation~(\ref{eqn:Prop2}) is equivalent to
 \begin{equation}\label{eqn:Prop2prim}
 a_v^{(j+1)}=a_v^{(j)}+a_{v-M_{\rho_j}}^{(j)}+
 a_{v-2M_{\rho_j}}^{(j)}+\ldots+
 a_{v-(\ell_j-1)M_{\rho_j}}^{(j)}\,.
 \end{equation}
 If $v\in S^{\K}_C$ cannot be represented as
 $M_{\rho_j}+b$ with $b\in S^{\K}_C$,
 one has $a_v^{j+1}=a_v^{j}$, what is compatible with~(\ref{eqn:Prop2prim}).
 Let $v=\ell M_{\rho_j}+b$ with
 $b\in S^{\K}_C$,
 $1\le \ell<\ell_j-1$,
 and $v$ cannot be represented as $(\ell+1) M_{\rho_j}+b'$ with $b'\in S^{\K}_C$.
In this case Equation~(\ref{eqn:Prop2prim}) reduces to
 \begin{equation}\label{eqn:Prop2primprime}
 a_v^{(j+1)}=a_v^{(j)}+a_{v-M_{\rho_j}}^{(j)}+
 a_{v-2M_{\rho_j}}^{(j)}+\ldots+
 a_{v-(\ell-1)M_{\rho_j}}^{(j)}+a_b^{(j)}\,
 \end{equation}
 (all other summands are equal to zero).
 Since, for $0\le\ell'\le\ell-1$,
 $v-\ell'M_{\rho_j}=M_{\rho_j}+b^*$ with
 $b^*\in S^{\K}_C$ and, due to
 Lemma~\ref{lemma:5_in_Th2} (Lemma \ref{lemma:1_in_Th2} for $j=1$),
 $a_{M_{\rho_j}}^{(j)}=[\K_{j-1}:\K_0]$,
 all the summands
 in~(\ref{eqn:Prop2primprime}) except the last one are equal to $[\K_{j-1}:\K_0]$.
 Thus in this case
 (\ref{eqn:Prop2prim}) follows from
 Lemma~\ref{lemma:4_in_Th2}.

 Let $v=(\ell_j-1) M_{\rho_j}+b$ with $b\in S^{\K}_C$.
 There are two possibilities.
 \begin{enumerate}
  \item[1)] The value $v$ can be represented as
  $M_{\varkappa}+b'$ with $b'\in S^{\K}_C$,
  $\varkappa>\rho_j$. In this case
  (\ref{eqn:Prop2prim}) follows from
  Lemmas~\ref{lemma:5_in_Th2}
  and~~\ref{lemma:6_in_Th2}.
  \item[2)] The value $v$ cannot be represented in the indicated form. In this case
 (\ref{eqn:Prop2prim}) follows from
 Lemma~\ref{lemma:4_in_Th2} again.
 \end{enumerate}
 \end{proof}

 \begin{remark}\label{rem:after_Prop2_in_Th2}
  One can see that the proof of Proposition~\ref{prop:Prop2_in_Th2}
  gives a little bit more strong statement
  which we shall use in Section~\ref{sect:PS-divis}.
  Let $A'_j$ be
the set of collections $\{k_{\sigma}\}$ such that $k_{\sigma}=0$
for all $\sigma$
on the geodesic from $\rho_j$ to $\delta_C$ {\bf excluding (!)} $\rho_j$ and
including $\delta_C$.
For $j=s+1$ (if $s<+\infty$), we assume this geodesic to be empty.
Let the series
$\widetilde{P}^{(j)}_C(t)=\sum_{v=0}^{\infty} \widetilde{a}^{(j)}_v t^v$ be defined by
$$
\widetilde{a}^{(j)}_v=\max_{\{k_{\delta}\}\in A'_j:\nu(\{k_{\sigma}\})=v}d^{\{k_{\sigma}\}}.
$$
Then one has
\begin{equation}
  \widetilde{P}^{(j)}_C(t)=P^{(j+1)}(t)=P^{(j)}(t)
  \cdot\frac{1-t^{\ell_j M_{\rho_j}}}{1-t^{M_{\rho_j}}}.
  \end{equation}
 \end{remark}

 \begin{Proof} {\bf of Theorem~\ref{theo:theo2}.}
 If $\rho_1\ne\sigma_0$, the statement follows directly
 from Propositions~\ref{prop:Prop1_in_Th2}
 and~\ref{prop:Prop2_in_Th2}.
 If $\rho_1=\sigma_0$,
 Lemmas~\ref{lemma:3_in_Th2}--\ref{lemma:6_in_Th2}
 and thus Proposition~\ref{prop:Prop2_in_Th2}
 hold for $j\ge 2$.
 Proposition~\ref{prop:Prop1_in_Th2}
 does not hold. However, it is possible to compute
 the series $P^{(2)}_C(t)$ directly
 using Lemmas~\ref{lemma:3_in_Th2}--\ref{lemma:6_in_Th2}.
 Thus the statement of Theorem~\ref{theo:theo2} holds in this case as well.
\end{Proof}

\begin{remark}
Any series $\zeta(t)= 1+ c_1 t + c_2 t^2 + \cdots\in \Z[[t]]$ can be in a unique way
represented as
$\zeta(t) = \prod_{m=1}^{\infty} (1-t^{m})^{s_m}$ with $s_m\in \Z$. We say that the
series is cyclotomic if this product/ratio of the binomials is finite.
There is a natural generalization of this notion for series in several variables.
Usually in previous computations (say, e.g., \cite{Duke}) the Poincar\'e series
appeared to be cyclotomic. This is not so for the Poincar\'e series
$P^{\K}_C(t)$ in case II, i.\,e., when $s=\infty$.
\end{remark}

\begin{remark} Let us assume that the curve valuation corresponds to Case I.
Let $c$ be the conductor of the
semigroup $S^{\K}_{C}$ (i.e., the minimal element of $S^{\K}_{C}$ such that
$c +\Z_{\ge 0}\subset S^\K_{C}$). It is known that $c$ can be expressed as
$$
c = \sum_{i=1}^g (N_i -1)M_{\sigma_i} - M_{\sigma_0} +1 \; .
$$
Let
$\ell = \ell_1 \cdots \ell_s = [\K_s : \K]$ and
$$
\Delta := c + \sum_{j=1}^s (\ell_j-1)M_{\rho_j} =
\sum_{i=1}^g (N_i -1)M_{\sigma_i} - M_{\sigma_0} +1
+ \sum_{j=1}^s (\ell_j-1)M_{\rho_j}\; .
$$
The computation of the Poincar\'e series gives that
$a_v = \dim J^{\K}(v)/J^{\K}(v+1) = \ell$ for all $v\ge \Delta$ and
$a_{\Delta-1} < \ell$. Moreover,
one has the following symmetry property:

\begin{corollary}
For all $v\in \Z$ one has
\begin{equation}\label{eq:symmetry}
\dim \left({J^{\K}(v)}\left/{J^{\K}(v+1)}\right.\right) + \dim\left(
{J^{\K}(\Delta-1-v)}\left/{J^{\K}(\Delta-v)}\right.\right) = \ell\; .
\end{equation}
\end{corollary}

Following the notations of Remark~\ref{rational}, the integral closure $\b{R}$
of the local ring $R = \K[[x,y]]/(F)$ of the algebroid curve defined by $F$
is a valuation ring and is
isomorphic to $\K'[[u]]$ ($u$ is an uniformizing parameter).
One has:

\begin{corollary}
The conductor ideal $\CC(\b{R}/R) = \{z\in \b{R} : z\b{R}\subset R\}$ is the ideal
$u^{\Delta}\b{R}$.
\end{corollary}

Note that the ring $R$ is Gorenstein (being plane) and the symmetry described by the equation
(\ref{eq:symmetry}) is just the characterization of the Gorenstein condition of
the ring $R$ (cf.~\cite{CDK}).

\end{remark}

\section{The Poincar\'e series of a divisorial valuation}\label{sect:PS-divis}  
Let $\nu = \nu_{\delta}$ be a divisorial valuation on $\OO_{\C^2,0}$.
This means that it is defined by a component
$E_\delta$ of the exceptional divisor of a modification of $(\C^2,0)$.
We assume that the corresponding $G$-invariant modification is
finite. Otherwise we are in a situation similar to Case III in the curve case
and the valuation $\nu_{\delta}$ is a multiple of
another one: corresponding to the component on which
the $G$-orbit of the blown-up point is infinite.
Let
$\pi: (W,\DD)\to (\C^2,0)$ be the minimal $G$-resolution of the valuation $\nu$. The
quotient of the dual graph of this resolution by the action of the Galois group $G$
looks like in figure \ref{fig5}.


\begin{figure}[h]
$$
\unitlength=0.50mm
\begin{picture}(120.00,40.00)(0,10)

\thinlines
\put(-30,30){\circle*{2}}
\put(-40,24){{\scriptsize ${\bf 1}=\sigma_0$}}

\put(-30,30){\line(1,0){20}}
\put(-8,30){\circle*{0.5}}
\put(-5,30){\circle*{0.5}}
\put(-2,30){\circle*{0.5}}
\put(1,30){\circle*{0.5}}
\put(4,30){\circle*{0.5}}
\put(5,30){\line(1,0){40}}

\put(-15,30){\line(0,-1){15}}
\put(-15,30){\circle*{2}}
\put(-15,15){\circle*{2}}
\put(-16.5,10){{\scriptsize$\sigma_1$}}
\put(-16,33){{\scriptsize$\tau_1$}}

\put(10,30){\line(0,-1){20}}
\put(10,30){\circle*{2}}
\put(10,10){\circle*{2}}
\put(11.5,7){{\scriptsize$\sigma_q$}}
\put(9,33){{\scriptsize$\tau_q$}}

\put(30,30){\circle*{2}}
\put(27,23){\scriptsize{$\rho_1$}}

\put(30,30){\line(1,1){10}}
\put(30,30){\line(1,2){8}}
\put(30,30){\line(0,1){10}}

\put(40,30){\line(0,-1){20}}
\put(40,30){\circle*{2}}
\put(40,10){\circle*{2}}
\put(41.5,7){{\scriptsize$\sigma_{q+1}$}}
\put(39,33){{\scriptsize$\tau_{q+1}$}}

\put(48,30){\circle*{0.5}}
\put(51,30){\circle*{0.5}}
\put(54,30){\circle*{0.5}}

\put(55,30){\line(1,0){40}}

\put(60,30){\line(0,-1){20}}
\put(60,30){\circle*{2}}
\put(60,10){\circle*{2}}

\put(75,30){\circle*{2}}
\put(72,23){\scriptsize{$\rho_j$}}

\put(75,30){\line(2,1){10}}
\put(75,30){\line(1,1){10}}
\put(75,30){\line(1,2){8}}
\put(75,30){\line(0,1){10}}

\put(90,30){\line(0,-1){20}}
\put(90,30){\circle*{2}}
\put(90,10){\circle*{2}}
\put(98,30){\circle*{0.5}}
\put(101,30){\circle*{0.5}}
\put(104,30){\circle*{0.5}}
\put(107,30){\circle*{0.5}}

\put(112,30){\line(1,0){20}}
\put(120,30){\circle*{2}}
\put(118,23){\scriptsize{$\rho_s$}}

\put(120,30){\line(2,1){10}}
\put(120,30){\line(1,1){10}}
\put(120,30){\line(0,1){10}}

\put(134,30){\circle*{0.5}}
\put(137,30){\circle*{0.5}}
\put(140,30){\circle*{0.5}}

\put(142,30){\line(1,0){8}}

\put(150,30){\circle{4}}
\put(150,30){\circle*{2}}

\put(155,28){{\scriptsize $\delta$}}

\end{picture}
$$
\caption{The resolution graph $\check \Gamma$ for a divisorial valuation.}
\label{fig5}
\end{figure}

This graph essentially coincides with the minimal resolution graph of a
$\K_{\delta}$-curvette at the component $E_{\delta}$ ($\K_{\delta}=\K_s$) with a
possible tail from $\tau_g$ to $\delta$ added. Here, as above,
$\sigma_i$, $i=0,1,\ldots, g$, are the dead ends of the graph $\check \Gamma$,
$\tau_i$, $i=1,\ldots, g$, are the rupture points and $\rho_j$, $j=1,\ldots,s$, are
the splitting points of the $G$-resolution.

\begin{remark}
Pay attention that $\rho_s\neq \delta$. One may define $\rho_{s+1}$ as $\delta$. In
this case the divisorial valuation $\nu_\delta$ coincides with the curve valuation
$\nu_{C^*}$ for $C^*$ being a curvette at the component $E_\delta$ at a non-algebraic point
of it. (Thus $\ell_{s+1}$ is assumed to be equal to $\infty$.)
\end{remark}

We shall keep the notations $\EE^{\{k_\sigma\}}_{\K,0}$,
$F^{\{k_\sigma\}}$, and $d^{\{k_\sigma\}}$ as before.
Let $C$ be a $\K_s$-curvette at the component $E_{\delta}$.

\begin{theorem}\label{theo:theo3}
One has
\begin{equation}\label{eq6}
P^{\K}_{\nu_{\delta}}(t) = P^S_C(t)\cdot
\frac{1}{1-t^{M_\delta}} \cdot
\prod_{i=1}^s\frac{1-t^{\ell_i M_{\rho_i}}}{1-t^{M_{\rho_i}}} \;.
\end{equation}
\end{theorem}

\begin{proof}
For the coefficient $a_v$ in the Poincar\'e series
$P^{\K}_{\nu_{\delta}}(t)=\sum a_v t^v$
one has
$a_v := \max d^{\{k_\sigma\}}$ where the maximum is taken over the collections
$\{k_\sigma\}$ such that $\nu(\{k_\sigma\}) = \sum k_\sigma M_\sigma=v$.
Let $A'$ be the set of collections $\{k_\sigma\}$, $\sigma\in \check \Gamma$ such
that $k_\delta=0$ and let
$P'(t)= \sum a'_{v}t^v$ with
$a'_v := \max d^{\{k_\sigma\}}$ where the maximum is taken over the collections
$\{k_\sigma\}\in A'$ such that $\nu(\{k_\sigma\}) = v$.

\begin{lemma}\label{lemma:1_in_Th3}
 One has
 \begin{equation*}
P'_{\nu}(t) = P_C^{\K}(t)=P^S_C(t)\cdot
\prod_{j=1}^s\frac{1-t^{\ell_j M_{\rho_j}}}{1-t^{M_{\rho_j}}}
\; .
\end{equation*}
\end{lemma}

\begin{proof}
The statement is a direct consequence
of Remark~\ref{rem:after_Prop2_in_Th2}.
\end{proof}

\begin{lemma}\label{lemma:2_in_Th3}
  Let $v\in S^{\K}_C$ be of the form
  $v=\ell M_{\delta}+b$ with $b\in S^{\K}_C$
  and 
such that $b-M_{\delta}\notin S_C^{\K}$.
  Then
  $a_{v}=\ell[\K_{s}:\K_0]+a'_{b}$.
 \end{lemma}

 \begin{proof} (Cf. proof of
 Lemma~\ref{lemma:4_in_Th2}.)
 Let us represent a collection $\{k_{\sigma}\}$
 with $\nu(\{k_{\sigma}\})=v$ in the form
 $\{k_{\sigma}\}=\{k^1_{\sigma}\}+\{k^2_{\sigma}\}$,
 where $k^1_\delta=\ell'$ ($\ell'\le\ell$),
 $k^1_\sigma=0$ for $\sigma\ne \delta$,
 $k^2_{\delta}=0$.
 A function $f\in\calE_{\K^2,0}^{\{k_{\sigma}\}}$
 can be represented as $f=f_1f_2$ with
 $f_r\in\calE_{\K^2,0}^{\{k^r_{\sigma}\}}$, $r=1,2$.
 Let $h_r$, $r=1,2$, be fixed functions from
 $\calE_{\K^2,0}^{\{k^r_{\sigma}\}}$, $h=h_1h_2$.
 As in Lemma~\ref{lemma:4_in_Th2},
 we shall use a coordinate on $E_{\delta}$
 such that the intersection point with the previous component from the geodesic
 to $\sigma_0$ is equal to $\infty$ and we shall take the function $h_1$ described there (that
 is, the $\ell'$th power of the left hand side of the equation of the $G_{s}$-curvette at the point $w=1$ of the component $E_{\delta}$).
 One has $\Psi={\widetilde f}/{\widetilde h}=
 \Psi_1\Psi_2$, where
 $\Psi_r={\widetilde f_r}/{\widetilde h_r}$.
 The function $\Psi_2$ is constant on the component $E_\delta$.
 The maximal possible dimension (over $\K$)
 of the set of its values on this component is $a'_{b'}$ ($b'=v-\ell'M_{\delta}$).
 The function $\Psi_1$ is constant on the union of all the components $E_\sigma$
for $\sigma$ from the
 geodesic between $\sigma_0$ and
 $\delta$, excluding the last one.
 Its restriction to the component $E_{\delta}$
 is equal to
 $$\frac{c \prod_{p=1}^{\ell'}(w-w_p)}{(w-1)^{\ell'}}$$ with $c\in\K$,
 where $w_p$ are the intersection points
 of the strict transform of the curve $\{f_1=0\}$
 with the component $E_{\delta}$. The set of these points
 is an arbitrary union of $G_s$-orbits.
 Therefore, up to the fixed multiple
 $(w-1)^{\ell'}$, the restriction of $\Psi$ to the component $E_{\delta}$
 is equal to the restriction of $\Psi_2$
 multiplied by $c\prod_{p=1}^{\ell'}(w-w_p)$.
 The dimension (over $\K$) of the space of functions of the latter form is
 $\ell'[\K_s:\K]$ (since the collection of the symmetric polynomials of the points
 $w_p$ is a generic point of $\K_s^{\ell'}$).
 Therefore the (maximal) dimension of the set of these functions on $E_{\delta}$ is
$a'_{b'}+\ell'[\K_s:\K_0]$.
 Thus it is maximal only for $\ell'=\ell$
 and is equal to $\ell[\K_{s}:\K_0]+a'_{b}$.
 \end{proof}

 Just in the same way as
 Proposition~\ref{prop:Prop2_in_Th2}
 we get the following statement.

 \begin{proposition}\label{prop:in_divis}
 One has
 $$
 P^{\K}_{\nu_{\delta}}(t)=P'(t)\cdot\frac{1}{1-t^{M_\delta}}
 =P^{\K}_{C}(t)\cdot\frac{1}{1-t^{M_\delta}}.
 $$
 \end{proposition}

Lemma~\ref{lemma:1_in_Th3} and Proposition~\ref{prop:in_divis}
imply the statement of Theorem~\ref{theo:theo3}.
\end{proof}

\begin{remark}
 One can see that the Poincar\'e
  series of a curve or of a divisorial valuation determines the combinatorial type
 of the minimal $G$-resolution of the
 valuation.
\end{remark}

\begin{remark}
The formula in Proposition \ref{prop:in_divis} relating the Poincar\'e series of the
divisorial valuation  $\nu$ and the Poincar\'e series of a $\K_s$-curvette $C$ at the
component $E_{\delta}$ obviously extends the known one for the complex case (see
\cite{Galindo-1995} or \cite{DG}). Even the formula itself and, to some extent, the
proofs have some technical resemblance with the relation of both series in the
case of several divisorial valuations in \cite{DGN}. However, the objects are rather
diferent, in \cite{DGN} all the considerations are in
the pure complex analytic  framework and the proofs use only properties of
pencils (another way to treat
meromorphic functions like in this paper) in the complex setting.

\end{remark}


Addresses:

A. Campillo and F. Delgado:
IMUVa (Instituto de Investigaci\'on en
Matem\'aticas), Universidad de Valladolid,
Paseo de Bel\'en, 7, 47011 Valladolid, Spain.
\newline E-mail: campillo\symbol{'100}agt.uva.es, fdelgado\symbol{'100}uva.es

S.M. Gusein-Zade:
Moscow State University, Faculty of Mathematics and Mechanics,
Moscow, Leninskie Gory 1, GSP-1, 119991, Russia.\\
\& National Research University ``Higher School of Economics'',
Usacheva street 6, Moscow, 119048, Russia.
\newline E-mail: sabir\symbol{'100}mccme.ru

\end{document}